\pdfoutput=1
\documentclass[12pt,a4paper]{article}

\usepackage[T1]{fontenc}
\usepackage[utf8]{inputenc}
\usepackage[margin=2.8cm]{geometry}
\usepackage{mathptmx}      
\usepackage{microtype}     
\usepackage{amsmath,amssymb,amsfonts,amsthm}
\usepackage{bm}
\theoremstyle{definition}

\usepackage{enumitem}

\usepackage{graphicx}
\usepackage{subcaption}
\usepackage{float} 

\usepackage{booktabs}
\usepackage{multirow}
\usepackage{makecell}
\usepackage{xcolor}

\setlist[itemize]{leftmargin=1.2cm}

\usepackage{cite}

\usepackage[hidelinks]{hyperref}

\begin{document}
\title{A Three Axis Evaluation Framework for Mapper Algorithms}

\author{
Annesha Sen\thanks{sen.annesha11@gmail.com}, Shivam Singh\thanks{singh.shivam18bsc083@gmail.com} and S.P. Tiwari\thanks{sptiwarimaths@gmail.com}\\
Department of Mathematics \& Computing\\
Indian Institute of Technology(ISM), Dhanbad-826004, India }
\date{}
\maketitle

\begin{abstract}
Mapper is a well-known tool in topological data analysis, which visualizes and summarizes high-dimensional data. However, its output is sensitive to choices of lens functions, cover parameters, and clustering strategies, making evaluation challenging. Most works that have attempted to evaluate the Mapper algorithm have done so visually. In this paper, we review a roadmap for assessing Mapper algorithms along three complementary axes: stability, cluster quality, and topological shape preservation. We analyze Mapper and its variants on synthetic datasets and the UCI Digits dataset. These modes include topological explosion at high resolutions. Our findings indicate that these axes of evaluation are often in tension and that no single Mapper variant performs optimally across all three. This review provides practical guidelines for choosing Mapper variants and identifies open challenges toward a principled Mapper analysis.
\end{abstract}


\section{Introduction}

Mapper~\cite{10,33} has become a widely used tool in topological data analysis (TDA), which is used mainly for visualizing high-dimensional datasets and has found many applications across biology, neuroscience, climate related issues ~\cite{76,77,75} materials science, and machine learning~\cite{32,33} and deep learning~\cite{71,16,19,15,39}. 
In biology and medicine, Mapper revealed biologically significant subgroups in high-dimensional data. For example, breast cancer subtypes that are not apparent with standard analyses \cite{45,46}. These ideas have been extended to diverse biological domains. In viral evolution, TDA has been used to reveal non-obvious evolutionary structure ~\cite{54,84,85}. In epidemiology, Mapper has been applied to capture global and temporal structure in influenza incidence data~\cite{53} and to study the spread of coronavirus~\cite{55}. In microbiome research, Mapper-based frameworks such as Two-Tier Mapper~\cite{21} and integrative TDA approaches~\cite{47,48,14,34} have been developed to preserve variance and highlight geographical or strain-level signatures in high-dimensional metagenomic data~\cite{49}. In ~\cite{57}, topological methods to single-cell RNA-seq to track differentiation trajectories and identify distinct developmental stages and molecular markers was applied.
Mapper has also been employed in neuroscience to find hidden structure in brain data~\cite{38}. Several variants, such as NeuMapper~\cite{83} (scalable, multiscale brain exploration), Temporal Mapper~\cite{87} (state-transition networks), and DeMapper~\cite{38}, have been developed to analyze brain network dynamics in clinical populations. Applications to ADHD include population-level disease stratification using subject-level graphs constructed from functional connectivity features~\cite{86} and subject-specific dynamical network analysis that constructs brain-state transition networks and quantifies topological organization via network centrality measures~\cite{87}. Together, these studies illustrate complementary uses of Mapper in neuroimaging, spanning population-level stratification and subject-specific dynamical analysis. For broader reviews of Mapper applications, see \cite{51} and the upcoming comprehensive review \cite{30}. 

Another line of work explores the use of the Mapper algorithm on neural network. An example is the Deep Graph Mapper~\cite{71}, which combines Mapper with Graph Neural Networks (GNNs) to perform graph pooling and visualization in a topologically grounded manner. Here GNN is employed as a learnable lens function that maps graph nodes into a low-dimensional space, over which the Mapper construction is applied to generate a coarsened graph representation. This highlights the potential of Mapper as a unifying topological framework for graph pooling and visualization within neural network models.

In materials science, Persistent Homology has been used to uncover hidden geometric structures in amorphous solids, revealing unique patterns in atomic configurations~\cite{74}. These atomic arrangements are converted into persistence diagrams~\cite{18}, which capture ring and cavity structures that are not detected using traditional machine learning methods (cf.~\cite{32,33}). In particular, TDA has been applied to financial systems to detect market crashes~\cite{69}. These findings suggest that TDA-based descriptors can serve as informative indicators of structural changes in complex dynamical systems. Persistent Homology has also been used to analyze brain dynamics, where preserving spatial and temporal structure is essential. Summarizing whole-brain activity graphs helps detect changes between different mental health states and captures individual differences in neural dynamics (cf.~\cite{68,83}). A systematic review of Persistent Homology and its sensitivity to parameter changes is presented in~\cite{33}. Recent advances in topological data analysis, including methodological developments and theoretical insights, are comprehensively discussed in~\cite{82,79,27,40,80}. Applications of TDA to real-world problems spanning multiple domains are presented in~\cite{78,44,73,36,50,37,56}.

Despite its widespread use, Mapper outputs are highly sensitive to choices of lens, cover, and clustering parameters. It is because the applied Mapper papers analyze the output mainly by visual intuition rather than quantitative criteria. It relies on the Homology of the resulting graph, which can be highly sensitive to parameter choices. Small perturbations in these choices can yield qualitatively different graph structures, leading to contradictory interpretations of the same dataset. A simple illustration of this issue is the choice of an inappropriate filter function, which can map the data into a space that is not aligned with the underlying structure of interest, leading to misleading Mapper representations. Thus, it becomes a challenge to reproducibility, particularly in scientific domains where Mapper-based summaries are used to support biological or clinical hypotheses. This review, therefore, focuses on developing a framework for understanding Mapper in the following three different perspectives, which will help in getting an idea for selecting parameters while doing Mapper analysis on datasets:

\begin{enumerate}[label=(\roman*)]
  \item \textbf{Stability}: Robustness of the Mapper graph to data sampling, and parameter perturbations.
  \item \textbf{Cluster quality}: Internal quality of clusters and their impact on graph interpretability.
  \item \textbf{Topological preservation/ Shape Preservation}: How faithfully different Mapper variants preserve the {true topology} and shape of the datasets.
\end{enumerate}
This review is organized around a three-axis evaluation framework for Mapper algorithms, focusing on stability, cluster quality, and topological or shape preservation. Section~2 introduces the necessary preliminaries, including formal definitions of Mapper, Reeb graphs, and related topological concepts, as well as the datasets and ground-truth constructions used throughout the study. 
Section~3 investigates Mapper from a stability perspective, combining theoretical insights with empirical and statistical stability experiments, and concludes with a discussion of open problems in stability analysis. Section~4 is devoted to cluster quality, where we review classical cluster validation theory, analyze how cluster quality arises in Mapper constructions, and present a series of experiments examining the effects of parameters, lens functions, and clustering choices. 
Section~5 studies topological and shape preservation, distinguishing Reeb-theoretic topology preservation from metric shape preservation, and compares different Mapper variants through controlled experiments on synthetic datasets with known ground truth. 
Section~6 provides an integrative comparison by analyzing the relationships between stability, cluster quality, and preservation, and offers practical insights and decision guidelines for selecting Mapper variants. 
Finally, Section~7 concludes the paper by summarizing key findings and outlining directions for future research. All findings in this work are obtained on the tested datasets, lens functions, and parameter ranges, and should be interpreted as systematic patterns within this scope rather than universal guarantees. 
\section{Preliminaries}
\label{sec:mapper_variants}
This section introduces the terminology and definitions used throughout the paper. It also lists the Mapper variants considered in this study. The section is divided into four subsections: Subsection \ref{2.1} gives a quick overview of basic definitions; Subsection \ref{2.2} introduces the Mapper variants used in this paper; Subsection~\ref{2.3} describes the datasets used in the paper along with their corresponding ground-truth values, which are employed consistently across all experiments; and, Subsection~\ref{2.4} validates these ground-truth definitions. 

\subsection{Background on Simplicial and Mapper Constructions}
\label{2.1}

This section provides the conceptual framework for the analysis in this paper. It highlights the definitions relevant to understanding the structure, interpretation, and limitations of Mapper outputs.

The beginning of topological data analysis lies in the notion of a simplex.
Given a collection of affinely independent points
$v_0, v_1, \dots, v_k \in \mathbb{R}^n$, their convex hull defines a
$k$-simplex, representing the simplest $k$-dimensional building block.
Simplices naturally generalize points, line segments, triangles, and
their higher-dimensional analogues, and they serve as the fundamental
units from which more complex topological objects are constructed (cf., ~\cite{42}).

A finite collection of simplices assembled in a consistent manner forms a simplicial complex. Specifically, a simplicial complex is required to contain all faces of its simplices and to ensure that intersections occur only along shared faces. This combinatorial structure provides a discrete
representation of a continuous space, while still retaining essential topological information such as connectivity and the presence of cycles. The dimension of a simplicial complex is defined as the maximum dimension among its constituent simplices (cf., ~\cite{42}).

Often, in order to analyze such complexes on different levels of detail, it is useful to consider their skeletons. The $k$-skeleton of a simplicial complex consists of all the simplices of dimension at most $k$. Of particular importance in practice is the $1$-skeleton, which includes only vertices and edges. The $1$-skeleton can be interpreted as a graph and forms the primary object of visualization in most Mapper-based analyses.

Simplicial complexes come with topological invariants quantifying features such as simplicial homology. For each dimension $k$, one builds a chain group generated by the $k$-simplices and defines a boundary operator that encodes how these simplices are to be glued together. The resulting homology groups measure the failure of boundaries to fully capture cycles, hence identifying topological structures such as connected components, loops, and higher-dimensional holes. The rank of the $k$th homology group, known as the $k$th Betti number, gives a compact encoding of all these features: $\beta_0$ counts connected components, $\beta_1$ counts independent loops, and $\beta_2$ counts enclosed voids (cf., ~\cite{42}).

Another important combinatorial construction is the nerve of a cover. Given an amount of coverage in a space, the nerve encodes the intersection pattern among the cover elements as a simplicial complex. Under suitable conditions, the nerve preserves the homotopy type of the underlying space, making it a central tool for translating continuous geometric information into a discrete topological representation (cf., ~\cite{42}).

The Mapper algorithm builds directly upon this idea. Given a finite metric space $(X,d)$, a user-chosen lens function $f: X \to \mathbb{R}^d$, and an open cover of the image $f(X)$, Mapper first pulls back each cover element to subsets of the data and then applies a clustering procedure within each subset. The resulting clusters serve as vertices of a simplicial complex, and simplices are introduced whenever clusters overlap. In practice, the $1$-skeleton of this complex, commonly referred to as the Mapper graph, is used for visualization and qualitative analysis.

Mapper is often interpreted as a discrete approximation of the Reeb graph \cite{28,66} associated with the same lens function. The Reeb graph is obtained by collapsing each connected component of a level set of a continuous function into a single point, thereby capturing how connectivity evolves along the function values. While Mapper converges to the Reeb graph under ideal conditions, finite sampling, clustering choices, and cover parameters introduce nontrivial distortions.

In this context, we distinguish between \emph{topological preservation}
and \emph{shape preservation}. Topological preservation refers to the
agreement of homological features such as Betti numbers between the Mapper output and the reference object (e.g., a Reeb graph or known ground truth), up to discretization effects. Shape preservation, on the other hand, emphasizes the stability of metric connectivity and neighborhood structure independent of the chosen lens function.

A recurring phenomenon in Mapper constructions is \emph{topological
explosion}, wherein increasingly fine covers lead to the rapid appearance
of spurious cycles in the Mapper graph, even when the underlying space has simple topology. This behavior highlights the delicate balance between resolution and robustness and motivates the need for principled parameter selection.

Finally, for synthetic datasets, the notion of \emph{ground truth
topology} plays a crucial role in evaluation. Since the generating space
and its homological invariants are known a priori, such datasets provide a
controlled setting for assessing whether Mapper faithfully recovers the
intended topological structure or introduces artifacts.

\subsection{Mapper Variants}
\label{2.2}
Since the introduction of the Mapper algorithm, a variety of extensions and alternative constructions have been proposed to adapt the method to different data characteristics and application domains. These variants modify specific components of the original Mapper pipeline, such as the choice of cover, clustering strategy, or graph construction, while preserving the overarching goal of producing a concise topological summary of complex data. In this section, we briefly introduce the Mapper variants considered in this review, focusing solely on their defining constructional features. Detailed analyses of their properties are deferred to subsequent sections. Several surveys have catalogued a growing family of Mapper variants distinguished by modifications to the cover, clustering strategy, or nerve construction (cf., ~\cite{30}). The following variants have drawn the attention of researchers:
\begin{enumerate}
    \item \textbf{Conventional Mapper}~\cite{10} constructs a simplicial complex as the nerve of clusters obtained from pullback sets induced by a cover of the lens space. Given a dataset, a user-defined filter function maps the data to a low-dimensional space, which is then covered by overlapping intervals. Each pullback set is clustered, and the resulting clusters form the vertices of the Mapper complex, with edges determined by non-empty intersections; the Mapper graph is obtained from its $1$-skeleton.
    
    \item \textbf{Ensemble Mapper}~\cite{11} denotes a class of approaches in which Mapper is applied multiple times using different parameter choices or data subsamples. It does not depend on a single Mapper output; instead, it combines all resulting graphs to produce a single output. This can be done in many ways, for example, via frequency-based summaries or consensus graphs, to highlight structural features that consistently appear across multiple Mapper graphs.
    
    \item \textbf{MultiScale Mapper}~\cite{25} extends the classical Mapper construction by including higher-order intersections among clusters in the nerve complex. Unlike the standard Mapper, which focuses on the $1$-skeleton, this adds higher-dimensional simplices arising from multiple overlaps. As a result, MultiScale Mapper provides a complete simplicial complex rather than only its $1$-skeleton.
    
    \item \textbf{Fuzzy Mapper ($F$-mapper)}~\cite{13} modifies the Mapper pipeline by replacing a hard, uniform cover with a cover defined via fuzzy memberships. It uses the concept of FCM~\cite{6} for selecting covers (cf.,~\cite{29,75}) ; see~\cite{59,41} for applications.
    
    \item \textbf{Ball Mapper}~\cite{12} creates a cover on the dataset rather than defining a filter function and pullback. It selects a set of centers and covers the data with balls of a specified radius. Nodes in the Ball Mapper graph correspond to these balls, and edges are added when the corresponding balls share data points. This construction emphasizes geometric proximity in the original data space; see~\cite{61,63,64} for applications.
\end{enumerate}

\subsection{Datasets and Ground Truth Definitions}
\label{2.3}

We consider three datasets throughout this paper: \emph{Swiss Roll with a Hole}, \emph{Noisy Circle}, and \emph{UCI Handwritten Digits}. All experiments were performed using Python~3.11 and NumPy~2.0.2. Mapper graphs were constructed using a custom Python implementation. KeplerMapper~2.1.0 was used only for reference and consistency checks, and not as the execution engine for the reported experiments. To avoid circular reasoning, we do not use clustering algorithms (e.g., DBSCAN)
to define ground truth. Instead, we rely on two independent references:
(i) the known parametric construction of synthetic manifolds commonly used in
topological data analysis~\cite{10}, and
(ii) Euclidean persistent homology computed directly on the embedded point
cloud~\cite{18,19}. A detail of the datasets is sketched below.

\begin{enumerate}
\item \textbf{Swiss Roll with Hole}.
The Swiss Roll (Figure~\ref{fig:swiss_noisy}) is generated via the parametric embedding $x = t\cos(t)$, $y = h$, $z = t\sin(t)$, where $t \in [1.5\pi, 4.5\pi]$ and $h \in [0, 21]$. A rectangular region is excised to introduce a puncture. The canonical dataset is generated from an initial pool of $N_{base}=2,500$ points; after hole excision, the final dataset contains approximately $2,100$ samples.

When analyzed as a point cloud in $\mathbb{R}^3$ using a Euclidean Rips filtration~\cite{19}, two robust $H_1$ features are observed: (i) an intrinsic loop corresponding to the puncture, and (ii) an extrinsic loop induced by the spiral embedding. To accommodate both interpretations, we distinguish:
\begin{itemize}
    \item \textbf{Intrinsic correctness:} $\beta_1 = 1$ (Puncture only),
    \item \textbf{Euclidean correctness:} $\beta_1 \in \{1,2\}$ (Puncture + Extrinsic Spiral).
\end{itemize}
Throughout our experiments, we adopt \textbf{Euclidean correctness}. Outputs with $\beta_1 > 2$ are classified as combinatorial noise.

\item \textbf{Noisy Circle}.
The canonical Noisy Circle (Figure~\ref{fig:swiss_noisy}) consists of $N=1,000$ samples drawn from a unit circle with Gaussian noise ($\sigma=0.05$). Although the raw geometry may contain multiple weak loops, Euclidean persistent homology consistently reveals a single dominant $H_1$ feature. We therefore treat $\beta_1 = 1$ as the effective Euclidean topological signal.

\begin{figure}[H]
    \centering
    \includegraphics[width=0.5\linewidth]{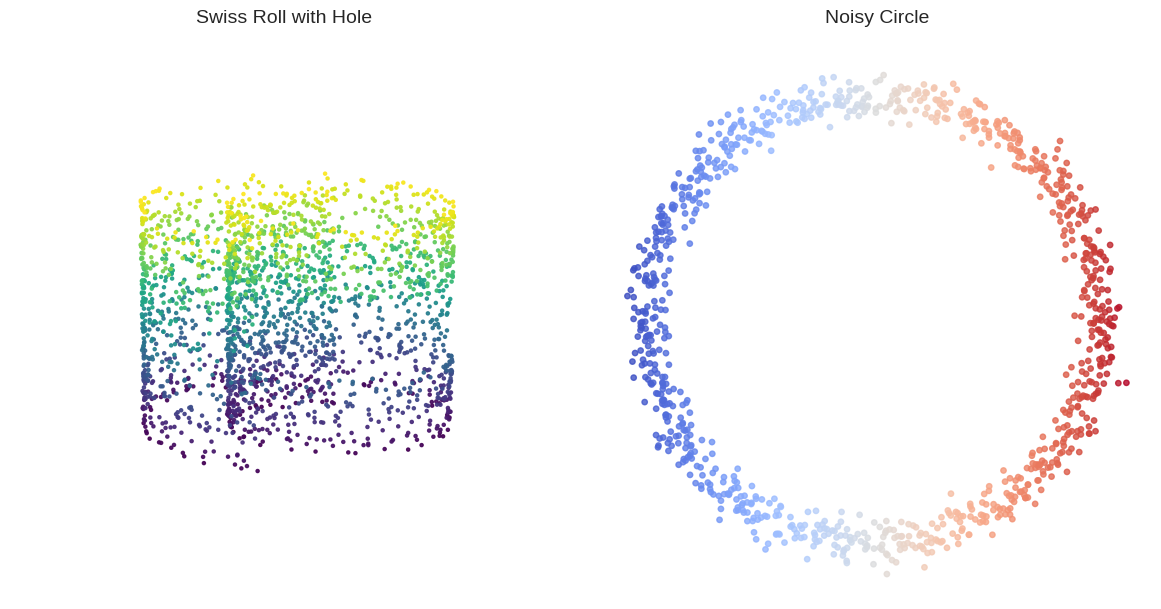}
    \caption{Swiss Roll with Hole and Noisy Circle Datasets.}
    \label{fig:swiss_noisy}
\end{figure}

\item \textbf{UCI Handwritten Digits}.
To probe stability in a high-dimensional regime, we use the 
\emph{Optical Recognition of Handwritten Digits} dataset from 
the UCI Machine Learning Repository, consisting of $N=1{,}797$ 
samples embedded in a $64$-dimensional feature space. The relatively 
small sample size compared to the ambient dimension makes this dataset 
well-suited for investigating density dilution and instability 
effects in Mapper constructions~\cite{16}. This dataset is used in Part~I (stability analysis) to document the high-dimensional regime 
behavior, and excluded from Parts~II--III due to the absence of 
ground-truth topology for comparative evaluation.

\end{enumerate}

\paragraph{Interpretation of $\beta_1$.}
Throughout this paper, the symbol $\beta_1$ is used in two distinct but complementary senses.
In Sections~\ref{2.3}--\ref{2.4}, $\beta_1$ denotes a \emph{topological ground-truth invariant}
computed via Euclidean persistent homology on the underlying point cloud.
In later sections, when analyzing Mapper outputs, $\beta_1$ refers to the
\emph{cycle rank of the Mapper graph 1-skeleton}, used as a diagnostic measure of
combinatorial complexity rather than a homology invariant.
These quantities are not expected to coincide in general, and their divergence
is a central object of study in this work.


\subsection{Validation of Ground Truth}
\label{2.4}

To verify that the ground-truth definitions are robust to sampling density and not artifacts of the full sample size, we performed a high-precision bootstrap stability analysis using Euclidean persistent homology (EPH) on random subsamples of the canonical data.

\paragraph{Methodology}
We performed \(T=20\) bootstrap trials with replacement. To rigorously test topological signal strength under density dilution, we subsampled the datasets as follows:
\begin{enumerate}
    \item \textbf{Swiss Roll:} Subsampled to \(N_{\mathrm{val}}=1{,}200\) (approx.\ 55\% of canonical size).
    \item \textbf{Noisy Circle:} Subsampled to \(N_{\mathrm{val}}=400\) (40\% of canonical size).
\end{enumerate}
Validating on these sparser subsets provides a stronger guarantee: if the topology is stable at reduced density, it is expected to be stable at the full canonical resolution used in Parts I--IV. We applied a strict persistence threshold, \(\tau=0.25\), to filter triangulation noise.

\paragraph{Findings}
Table~\ref{tab:gt_validation} summarizes the results.

\begin{itemize}
  \item \textbf{Noisy Circle:} Demonstrated complete stability (100\%), reliably retrieving \(\beta_1=1\).
  \item \textbf{Swiss Roll:} Attained full stability (100\%) within the Euclidean correctness interval \(\beta_1\in\{1,2\}\). This result indicates that \(N_{\mathrm{val}}=1{,}200\) together with the filter \(\tau=0.25\) effectively suppress combinatorial noise.
\end{itemize}

\begin{table}[htbp]
\centering
\caption{Validation of ground-truth topology using the subsampling method (\(T=20\), \(\tau=0.25\)). Stability at these low sampling densities indicates that the full canonical datasets are stable.}
\label{tab:gt_validation}
\begin{tabular}{lccc c}
\toprule
\textbf{Dataset} & \textbf{Validation \(N\)} & \textbf{Canonical \(N\)} & \(\boldsymbol{\beta_1}\) (observed) & \textbf{Success Rate} \\
\midrule
Noisy Circle & 400  & 1{,}000      & \([1,\,1]\) & 100\% \\
Swiss Roll   & 1{,}200 & \(\sim 2{,}100\) & \([2,\,2]\) & 100\% \\
\bottomrule
\end{tabular}
\end{table}

\subsubsection*{Persistence filtering and noise sensitivity}
\begin{figure}[htbp]
  \centering
  \begin{subfigure}[b]{0.48\textwidth}
    \centering
    \includegraphics[width=\textwidth]{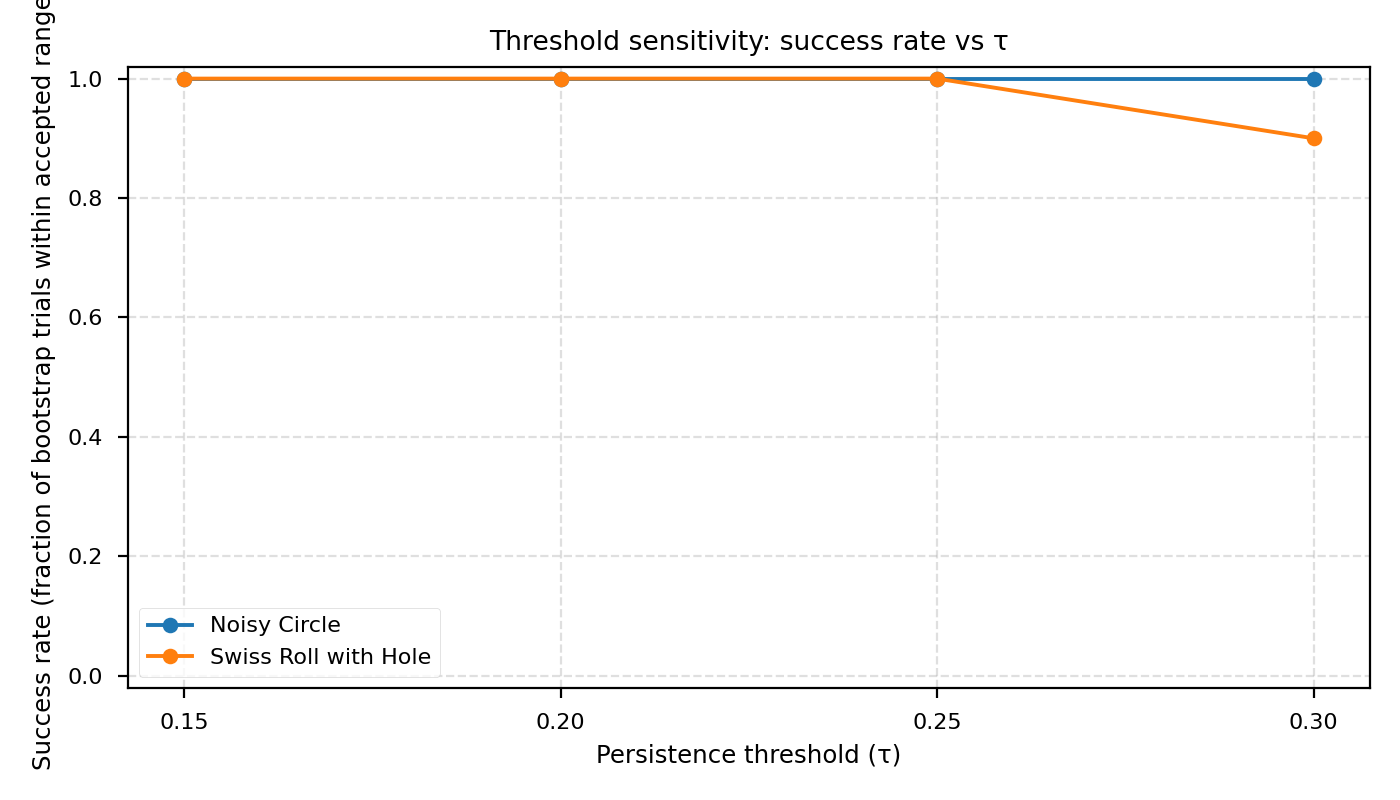}
    \caption{Success rate vs.\ \(\tau\).}
  \end{subfigure}
  \hfill
  \begin{subfigure}[b]{0.48\textwidth}
    \centering
    \includegraphics[width=\textwidth]{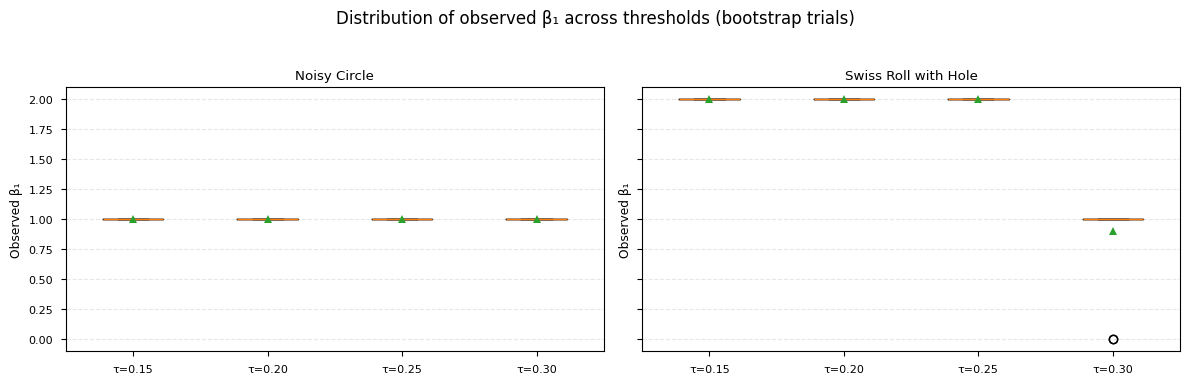}
    \caption{Distribution of observed \(\beta_1\) across resamples.}
  \end{subfigure}
  \caption{Threshold sensitivity analysis: (left) fraction of bootstrap trials whose observed \(\beta_1\) falls within the accepted range; (right) bootstrap distributions of \(\beta_1\) for each threshold.}
  \label{fig:tau_sensitivity}
\end{figure}

To remove short-lived homology classes caused by limited sampling and added noise, we filter the \(H_1\) barcodes using a persistence threshold \(\tau\) and only consider features with persistence \(>\tau\). To validate this parameter selection we performed a bootstrap sensitivity analysis with 20 resamples on the standard datasets, adjusting \(\tau\in\{0.15,0.20,0.25,0.30\}\). Success rates are reported in Table~\ref{tab:tau_sensitivity}, and a visual summary is shown in Figure~\ref{fig:tau_sensitivity}.

\begin{table}[htbp]
\centering
\caption{Sensitivity of persistence threshold to success rate. Success rate is the proportion of bootstrap trials for which the measured \(\beta_1\) lies within the tolerated values. For the Swiss Roll with Hole dataset we accept \(\beta_1\in\{1,2\}\); for the Noisy Circle we accept \(\beta_1=1\).}
\label{tab:tau_sensitivity}
\begin{tabular}{lcccc}
\toprule
Dataset & \(\tau=0.15\) & \(\tau=0.20\) & \(\tau=0.25\) & \(\tau=0.30\) \\
\midrule
Noisy Circle & 100\% & 100\% & 100\% & 100\% \\
Swiss Roll with Hole & 100\% & 100\% & 100\% & 90\% \\
\bottomrule
\end{tabular}
\end{table}

\begin{table}[htbp]
\centering
\caption{Evaluation of Mapper variants across multiple criteria. Checkmarks denote comprehensive evaluation; specialized symbols indicate framework-specific assessment or postponed analysis.}
\label{tab:variant_evaluation_axes}
\resizebox{\textwidth}{!}{%
\begin{tabular}{lcccc}
\toprule
\textbf{Variant} & \textbf{Stability Analysis} & \textbf{Cluster Quality} & \textbf{Topological Preservation} & \textbf{Primary Innovation} \\
\midrule
Conventional   & \(\checkmark\) & \(\checkmark\) & \(\checkmark\) & Reference baseline \\
\(F\)-mapper  & \(\checkmark\) & \(\checkmark\) & \(\checkmark\) & Adaptive cover \\
Ensemble       & \(\checkmark^{\dagger}\) & \(\checkmark\) & \(\checkmark^{\dagger}\) & Consensus aggregation \\
Ball Mapper    & \(\checkmark^{\ddagger}\) & -- & \(\checkmark\) & Geometric \(\varepsilon\)-net \\
Multiscale     & \(\checkmark^{\S}\) & -- & \(\checkmark^{\S}\) & Multi-resolution tower \\
\bottomrule
\end{tabular}%
}
\vspace{0.6em}

\footnotesize
\textsuperscript{\(\dagger\)}\textbf{Ensemble:}
Due to its meta-algorithmic nature (consensus aggregation across
multiple base-Mapper realizations), stability and topological analysis
is deferred to Part~IV (Section~\ref{sec:integrative}), where comparative
ensemble performance is evaluated against optimized
single-method baselines from Parts~I--III.

\textsuperscript{\(\ddagger\)}\textbf{Ball Mapper:}
Stability assessed via \emph{geometric robustness}
under Hausdorff-bounded noise (Theorem~3.5 in~\cite{12})
instead of parameter sensitivity.

\textsuperscript{\(\S\)}\textbf{Multiscale:}
Stability and topological preservation examined through an \emph{algebraic framework}
(persistence-module \(\varepsilon\)-interleaving; Theorem~4.8 in~\cite{25}) rather than empirical benchmarking.

\smallskip

\noindent\textbf{--} denotes an evaluation axis not relevant to the method's principal design goals.
\end{table}

Alongside the single-configuration Conventional and \(F\)-mapper baselines, we incorporate Ensemble Mapper as a meta-algorithmic parameter-selection approach, following Kang and Lim (2021)~\cite{11}, to evaluate whether consensus across multiple Mapper instances enhances cluster quality. Although Ball Mapper was first introduced as a preprint~\cite{12}, it has since been adopted in a range of applied TDA settings, including financial data analysis~\cite{63}, neural network diagnostics, and general-purpose data-science workflows supported by software implementations in Python, R, and Stata~\cite{61,64}. In this review, we therefore examine Ball Mapper as an algorithmic object and analyze its empirical stability, clustering behavior, and shape preservation independently of its publication status.

\section*{
\begin{center}
PART I
\end{center}
}
\section{Stability Perspective}
\label{Part I - Stability Perspective}

The first part addresses the notion of stability of Mapper graphs under perturbations.
It analyzes the stability of Mapper and its variants with respect to different types of
perturbations. This part is divided into four subsections:
Section~\ref{sec:mapper_stability_definitions} introduces the notion of stability from
multiple perspectives;
Section~\ref{sec:theoretical_mapper_stability} reviews the established theoretical
results for the Conventional Mapper and its variants;
Section~\ref{sec:stability_experiments} presents stability experiments conducted on the
datasets listed in Table~\ref{tab:variant_evaluation_axes};
and Section~\ref{sec:stability-open-problems} discusses open problems specific to
stability analysis.


\subsection{Stability Analysis}
\label{sec:mapper_stability_definitions}

The notion of stability in Mapper does not admit a single universal definition.
Depending on the object of comparison and the source of perturbation, stability
has been formalized and interpreted in several distinct ways. In this review,
we distinguish the following non-interchangeable notions.

\textit{(Note: The specific terms ``empirical stability'' and ``statistical
stability'' as distinct categories for Mapper analysis do not appear to be
established terminology in the Mapper/TDA literature, but we have introduced
these terms to make the discussion more reader-friendly.)}

\begin{itemize}
    \item \textbf{Empirical Stability} refers to the robustness of Mapper graph
    constructions under perturbations of parameters. In this setting, stability
    is assessed by systematically varying the parameters of the Mapper pipeline
    and observing the resulting changes in the graph structure. The main question
    here is: \emph{How sensitive is the Mapper output to small changes in its
    user-defined parameters?}
    
    The perturbation is done by varying cover resolutions, clustering algorithm
    parameters, etc. Stability is then quantified by comparing Mapper graphs
    obtained under these perturbations using graph-level measures, such as node
    and edge consistency, connected component preservation, or graph similarity
    metrics.

    \item \textbf{Statistical Stability} refers to the variability of Mapper
    outputs under variation in data samples. In contrast to empirical stability,
    the algorithmic parameters of Mapper are typically held fixed, and stability
    is assessed by perturbing the data itself. The main question here is:
    \emph{How does the Mapper graph change if one slightly changes the dataset?}
    
    The primary sources of variability here are sampling effects and measurement
    noise. Common approaches include bootstrap resampling and subsampling, where
    multiple datasets are generated from the same underlying distribution and the
    Mapper algorithm is applied independently to each sample. Stability is then
    measured by the consistency of the outputs.

   \item \textbf{Algorithmic Volatility vs. Systematic Bias (Diagnostic Distinction).} Besides determining the origin of perturbation, it is also important to distinguish the type of failure mode. This distinction does not introduce a new notion of stability but rather it describes the qualitative type of instability that arises depending on whether empirical or statistical perturbations occur. Volatility denotes unpredictable changes of topological invariants, such as randomly fluctuating Betti numbers following parameter changes. In contrast, bias denotes a consistent but incorrect model, such as a Mapper graph consistently suggesting the presence of a hole which in fact does not exist. As we will see in Section~\ref{sec:stability_ensemble}, distinguishing an ``unstable'' graph from a ``stably wrong'' one is an important point when testing ensemble methods.
   
\end{itemize}

\subsection{Theoretical Stability of Mapper Variants}
\label{sec:theoretical_mapper_stability}
In this subsection, we discuss the established theories of listed Mapper Variants.

\begin{enumerate}
    
\item  \textbf{Mapper Stability.}
The stability of the Mapper algorithm is of utmost importance from the data analysis perspective, since small perturbations in input data, sampling, or parameter choices can lead to substantially different Mapper graphs. This sensitivity arises from filter functions, cover constructions, and clustering procedures, and has motivated a diverse body of work aimed at formalizing and measuring instability.

A theoretical treatment is provided in~\cite{23}, where the first intrinsic numerical measure of instability for Mapper-type algorithms is developed. Mapper outputs were then represented as Mapper functions, enabling pointwise comparison of graphs obtained from independent samples. The instability of the Mapper algorithm on samples of size $n$ is defined as
\[
\mathrm{InStab}_{\mathrm{Mapper}}
\;:=\;
\mathbb{E}\!\left[
I\bigl(\{Q^{(i)}_{n_i}\}_{i=1}^t\bigr)
\right],
\]
where the expectation is taken with respect to the product probability measure
induced by $P$ (Probability Distribution) on pairs of independent samples in $X^n \times X^n$.

Complementary stability results have been developed from a topological perspective by relating Mapper to Reeb graphs and persistent homology. In~\cite{26}, Carrière and Oudot, as well as in~\cite{25}, Dey, Mémoli, and Wang studied the stability after connecting it to Reeb Graphs and showed that under certain conditions, a Mapper graph can be approximated by a discretized or pixelized Reeb Graph. Furthermore, the extended Bottleneck distance for analyzing stability is studied. In Multiscale Mapper~\cite{25}, the cover is modified in such a manner that the resulting Mapper output is a simplicial complex rather than a 1-skeleton, and then stability is analyzed using the Bottleneck distance.

Extending the above ideas to point cloud data, in~\cite{24}, Bungula and Darcy analyze Mapper stability through bifiltrations over cover resolution and clustering parameters. The work focuses on the clustering algorithm using DBSCAN and multi-parameter filtrations. Further, filtrations of Mapper graphs via clustering algorithms, free-border points as the main obstruction, and instability of one-parameter Mapper persistence were discussed. Using interleavings of persistence modules, it is shown that although one-parameter Mapper filtrations may be unstable, stability can be recovered by jointly varying cover size and clustering scale, yielding provable bounds under data perturbations.
\end{enumerate}

Several works~\cite{23,29} investigate stability through graph invariants and resampling-based approaches. Studies on the stability of Mapper graph invariants analyze how structural quantities, such as node counts, connected components, and loop statistics, behave under bootstrap sampling, providing empirical evidence of convergence as the sample size increases. These approaches do not provide theoretical guarantees, but they do offer practical insights into robustness and reproducibility.

Another theory, namely the Two-Tier Mapper (TTMap), was developed in~\cite{21}, which focuses on instability by modifying the Mapper construction itself. Specifically, it introduces a stability-oriented algorithm for high-dimensional biological data. By separating deviation modeling from topological clustering and employing a two-tier cover that captures both global and local structure, TTMap achieves robustness to normalization choices, sampling variation, and outliers. Importantly, TTMap is among the few Mapper-based methods to provide formal stability guarantees with respect to filter perturbations, domain variations, and sampling effects.

In~\cite{11}, Kang and Lim proposed an Ensemble Mapper framework in which multiple Mapper graphs are constructed over varying resolution parameters. Fitzpatrick et al.~\cite{22} further generalized this idea by treating Mapper as a family of graph constructions over a large parameter space and identifying stable regions via graph-level similarity measures. These ensemble approaches interpret stability as structural recurrence across parameter perturbations and subsampling, offering practical robustness at finite sample sizes, though without theoretical
convergence guarantees.

All the cited works reveal aspects of the stability of the conventional Mapper algorithm. Among these works, \cite{11,22,23,24} mainly focused on statistical stability, and \cite{21,25,26} primarily focused on empirical stability. Despite so many theories, a unified framework that optimizes lens function, cover, and a clustering algorithm for stability under perturbation remains lacking.

\begin{enumerate}
    
\item \textbf{MultiScale Mapper}:
As mentioned above, Multiscale Mapper refines the classical construction by retaining connected components of higher-order intersections. Existing results establish convergence and Reeb-theoretic fidelity under appropriate regularity conditions. However, a comprehensive stability theory comparable to that for conventional Mappers, including explicit bounds under parameter and data perturbation, remains a gap.

\item \textbf{Ensemble Mapper}:
Ensemble Mapper aggregates multiple Mapper outputs constructed under varying parameters or subsamples. It is stable under parameter and data perturbation, but the open areas here are the choice of filter function and the cover to be chosen.

\item \textbf{$F$-mapper}:
$F$-mapper replaces the cover as mentioned above. The whole stability analysis is an open area (cf., ~\cite{13}).

\item \textbf{Ball Mapper}:
Ball Mapper mainly focuses on shape preservation. However, it enjoys strong geometric stability guarantees. Dłotko~\cite{12} proved that the nerve of a union of balls is stable under Hausdorff perturbations (Theorem~3.5). Specifically, if two point clouds are $\delta$-close in Hausdorff distance, the topological features of their $\epsilon$-ball nerves are interleaved between scales $\epsilon$
and $3\epsilon + 2\delta$.

\end{enumerate}

\subsection{Stability Experiments}\label{sec:stability_experiments}
This section presents controlled \emph{illustrative experiments} designed to visualize the theoretical stability phenomena discussed in Section \ref{sec:mapper_stability_definitions}. Our objective is not to benchmark Mapper variants for performance or optimality, but to demonstrate how different notions of stability empirical and statistical manifest in finite sample regimes under canonical algorithmic choices, consistent with prior stability analyses of Mapper-type constructions \cite{8,23,26}.

All methods are implemented according to their original definitions \cite{10,12,13}. Reported numerical values are \emph{parameter and dataset dependent} and are interpreted qualitatively. 
Throughout Part I, $\beta_1$ refers to the cycle rank of the Mapper 1-skeleton graph, not simplicial homology of the full nerve. This measure captures combinatorial instability rather than intrinsic topological invariants.

\begin{enumerate}
   \item \textbf{Geometric Stability: Hausdorff Perturbations}
\label{sec:geometric_exp}

We first illustrate the geometric behavior of Ball Mapper under bounded perturbations (Section~3.1). We apply Ball Mapper to the Swiss Roll dataset ($N=2{,}500$) under varying Gaussian noise levels ($\sigma \in \{0.01, 0.05\}$), following the original construction of Ball Mapper \cite{12}. As shown in Figure~\ref{fig:geometric_heatmap}, the method exhibits the \emph{Lattice Effect}.

\begin{figure}[H]
    \centering
    \includegraphics[width=0.5\textwidth]{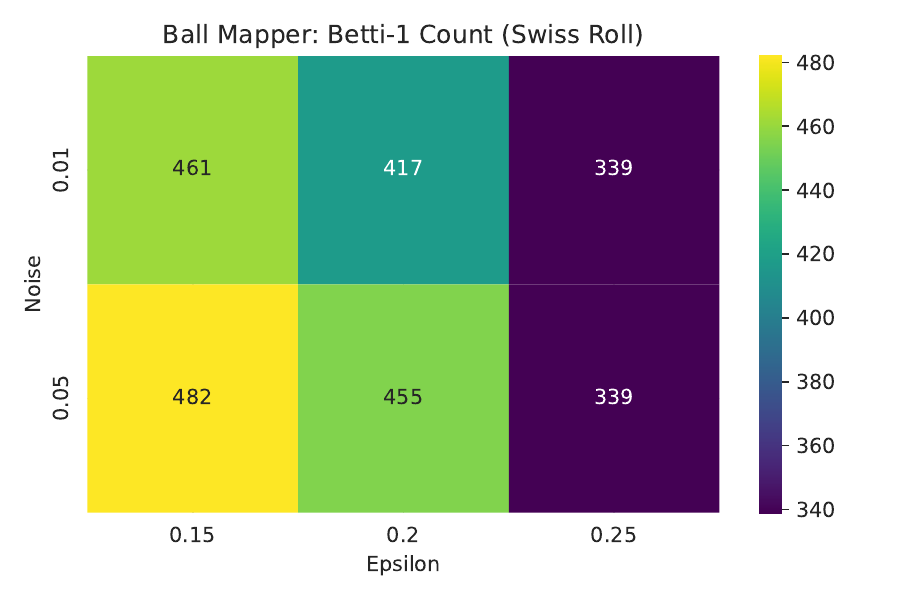}
    \caption{\textbf{Geometric Stability Heatmap (Ball Mapper).} 
    The heatmap displays the $\beta_1$ count for Ball Mapper on the Swiss Roll across varying noise levels ($\sigma$) and radius parameters ($\varepsilon$). The consistently high Betti numbers (ranging from 339 to 482) indicate that the nerve captures the geometry of the cover (the dense $\varepsilon$-net lattice) rather than the underlying manifold topology ($\beta_1=2$). This phenomenon, which we term the \emph{Lattice Effect}, reflects the packing density of balls on the manifold surface. Data source: \texttt{part1\_geometric\_stability\_canon.csv}.}
    \label{fig:geometric_heatmap}
\end{figure}

\paragraph{Key Observation.}
On the Swiss Roll manifold, Ball Mapper always generates around hundreds of 1-cycles. This is because of the packing density of balls on the manifold surface. It is geometrically consistent as it maintains the lattice at any amount of noise. It is topologically different from the manifold.

\item \textbf{Empirical Stability: Parameter Sensitivity}
\label{sec:empirical_exp}

For a clear understanding of the sensitivity of the algorithms towards parameters, we vary the resolution for Conventional Mapper and the scale parameter for Ball Mapper, following standard Mapper parameterizations \cite{8,10,26}.

\begin{figure}[htbp]
    \centering
    \includegraphics[width=0.95\textwidth]{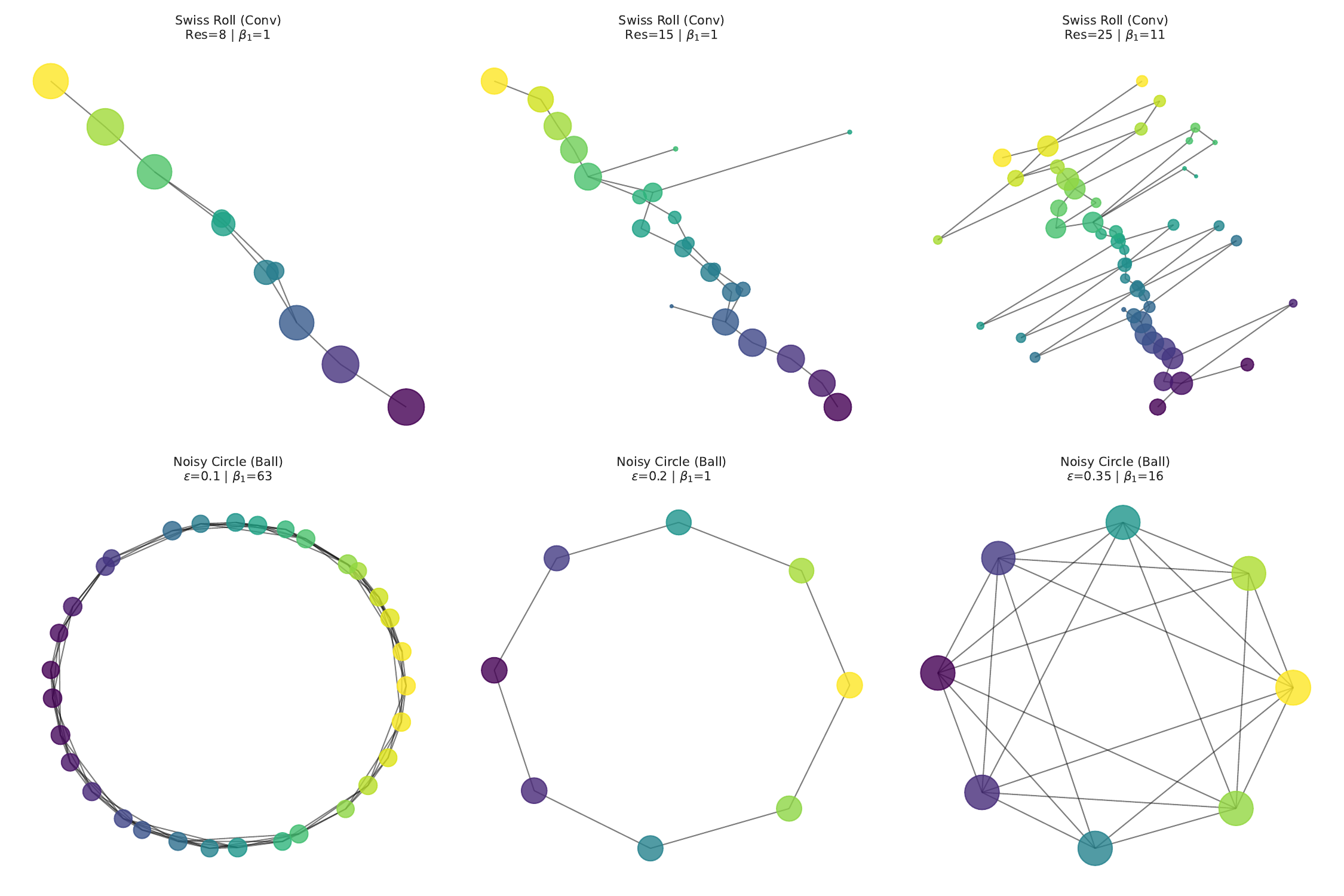}
    \caption{\textbf{Visual depiction of instability scenarios.} 
    (Top Row) \textbf{Conventional Mapper on Swiss Roll}: For a resolution of $\text{Res}=8$ and $15$, the graph captures the topology correctly ($\beta_1=1$). For a resolution of $25$, there is a topological explosion ($\beta_1=11$). 
    (Bottom Row) \textbf{Ball Mapper on Noisy Circle}: With $\varepsilon=0.1$, a dense lattice forms ($\beta_1=63$). At $\varepsilon=0.2$, the lattice collapses to correct topology ($\beta_1=1$). At $\varepsilon=0.35$, spurious cycles reappear ($\beta_1=16$).}
    \label{fig:instability_vis}
\end{figure}

\subsubsection*{High-Dimensional Behavior (UCI Digits)}
\label{sec:high_dim_exp}
We use the UCI Digits dataset ($d=64$, $N=1{,}797$) to illustrate how dimensionality interacts with cover construction\cite{12,29,70}.

\begin{table}[htbp]
\centering
\caption{\textbf{Illustrative high-dimensional behavior.}
Filter-based approaches shatter under projection, while Ball Mapper exhibits the \emph{Lattice Effect} in extreme high-dimensional form: the $\varepsilon$-balls form a connected structure but involve thousands of empty space, a feature compatible with ambient sparsity in high-dimensional metric spaces. Data source: \texttt{part1\_high\_dim\_canon.csv}.}
\label{tab:high_dim_failure}
\begin{tabular}{lccl}
\toprule
\textbf{Method} & \textbf{$\beta_0$} & \textbf{$\beta_1$} & \textbf{Qualitative Behavior} \\
\midrule
Classical Mapper & $12$ & $1$ & Projection Fragmentation \\
$F$-mapper & $12$ & $1$ & Smoothed fragmentation \\
Ball Mapper & $1$ & $11{,}054$ & Extreme Lattice Effect (high-dimensional) \\
\bottomrule
\end{tabular}
\end{table}

\paragraph{Analysis.}
Filter-based approaches (Conventional, $F$-mapper) rely on low-dimensional projection. Given this parameter configuration, projection groups data into distinct clusters with limited iterations ($\beta_0=12$) compared with Ball Mapper. Ball Mapper runs in the ambient space. Owing to high-dimensional space, $\varepsilon$-balls comprise a single connected component ($\beta_0=1$) with thousands of empty pockets ($\beta_1=11{,}054$) because of the concentration of measure principle of high-dimensional spheres.

\item \textbf{Statistical Stability: Bootstrap Resampling}
\label{sec:statistical_exp}

In order to assess the variability of statistics, we use the $T=50$ bootstrap samples on the Swiss roll dataset, following standard resampling-based stability analyses in TDA \cite{8,29}. The Coefficient of Variation (CV) is provided for the value of $\beta_1$, and the lower the percentage, the more stable it is than the mean.

\begin{figure}[htbp]
    \centering
    \includegraphics[width=0.5\textwidth]{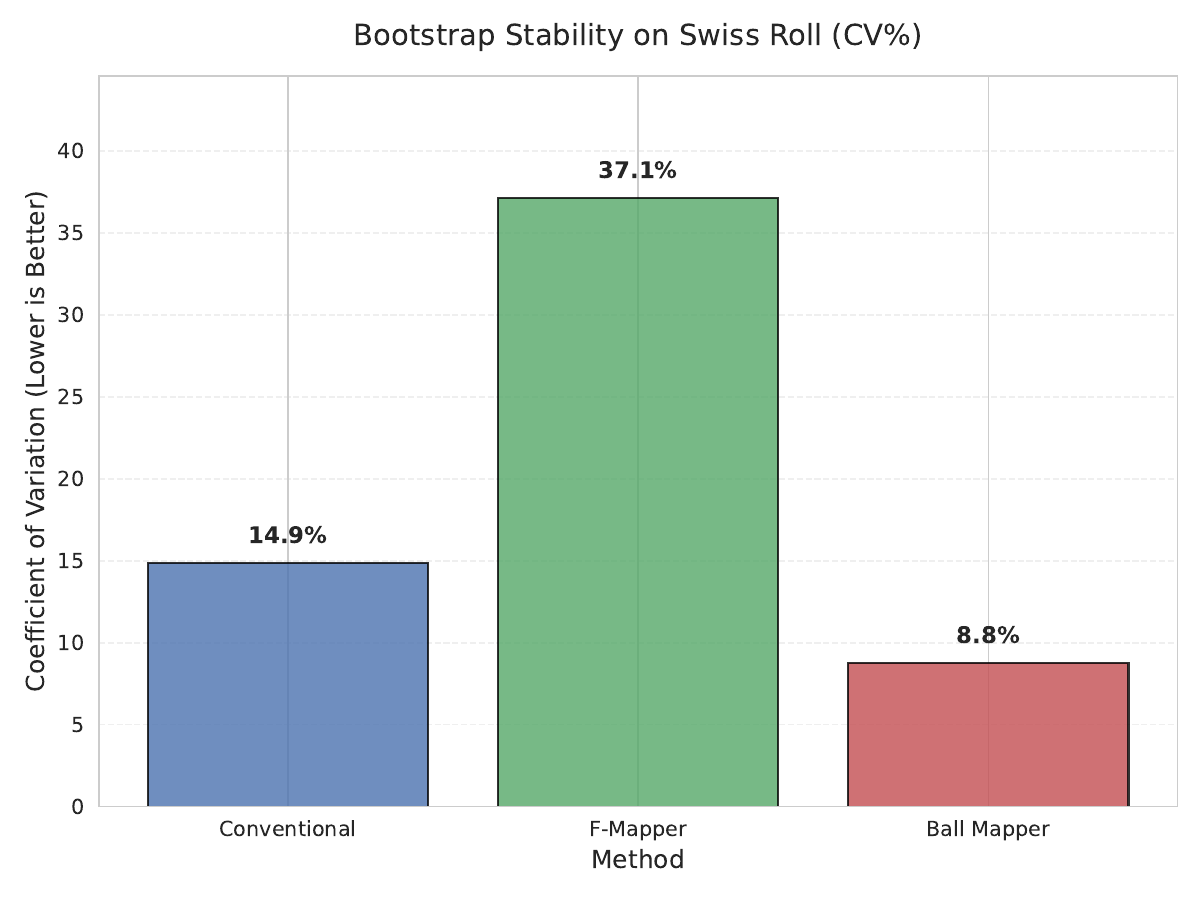}
    \caption{\textbf{Bootstrap variability of $\beta_1$ (Swiss Roll).} 
    \textbf{Ball Mapper} ($CV = 8.8\%$) is the most stable, as the lattice artifact is robust to resampling. 
    \textbf{Conventional Mapper} ($CV = 14.9\%$) shows low instability. 
    \textbf{$F$-mapper} ($CV = 37.1\%$) exhibits higher relative variance; this occurs because fuzzy memberships tend to suppress cycles \cite{13,41} (mean $\beta_1 = 3.6$, range $1$--$6$), making occasional detections statistically volatile relative to the low baseline. 
    Data source: \texttt{part1\_bootstrap\_canon.csv}.}
    \label{fig:bootstrap_stability}
\end{figure}
\end{enumerate}

\paragraph{Conclusion.}
These results highlight a dichotomy in stability as:
\begin{enumerate}
    \item \textbf{Metric Stability:} Ball Mapper is highly stable ($CV=8.8\%$) but stabilizes around a geometric artifact (the Lattice Effect) rather than the Reeb graph \cite{12,26}. While this is a limitation of the Reeb graph approximation, it provides geometric completeness for metric-based tasks, as we demonstrate in Part III (Section~5.3.4).
    \item \textbf{Filter Instability:} $F$-mapper's higher CV ($37.1\%$) illustrates the difficulty of tuning fuzzy membership parameters\cite{13,41}; the method vacillates between suppressing features entirely (bias) and detecting them, leading to statistical variance. Conventional Mapper provides a balanced trade-off ($CV=14.9\%$).
\end{enumerate}

\subsection{Open Problems in Stability Analysis}
\label{sec:stability-open-problems}

Despite substantial theoretical and empirical progress, stability in Mapper-based
constructions remains only partially understood. Several fundamental open problems
persist, both at the theoretical and practical levels.

\begin{itemize}
   \item      Existing stability results are fragmented across different notions, including
    empirical stability (parameter sensitivity), statistical stability (sampling
    variability). Therefore, a unified framework
    that simultaneously accounts for all these is still lacking.

    \item
    Most papers only analyze individual components of Mapper, such as the
    cover construction or the nerve approximation to Reeb graphs. But a combined
    effect of lens choice, cover resolution, overlap, and clustering parameters on the stability of the final Mapper graph has not been formally characterized.

   \item Stability theory of Conventional Mapper is well established, but that of other variants still lacks.

    \item 
    Current Mapper pipelines rely heavily on user-defined parameters. Developing principled, data-driven strategies for parameter selection that guarantee stability without extensive manual tuning is a significant open problem.

    \item 
    There is no universally accepted metric for comparing Mapper graphs across
    perturbations. Designing stability measures that are both computationally feasible
    and topologically meaningful remains an unresolved issue.

\end{itemize}

Collectively, these open problems highlight the need for stability analyses that move
beyond isolated components and toward end-to-end, data-adaptive, and theoretically
grounded Mapper constructions.

\section*{
\begin{center}
PART II
\end{center}
}

\section{Cluster Quality Perspective}
\label{Part II - Cluster Quality Perspective}

Cluster quality in Mapper refers to the quality of clusters that are obtained in the Mapper Graph. It plays an important role in determining the interpretability and reliability of Mapper outputs~\cite{10}. Because Mapper represents clusters of pullback sets as nodes in a graph, any degradation in clustering quality is immediately reflected in the resulting topological summary. This section is divided into three subsections: Section~\ref{cluster-theory} describes the cluster validation theory, Section~\ref{cluster-quality-in-mapper} describes cluster quality in Mapper and established theories in this area, Section~\ref{sec:cluster_quality_experiments} demonstrates related experiments, and Section~\ref{open-problems} discusses open problems.

\subsection{Cluster Theory}
\label{cluster-theory}
Cluster validation theory~\cite{2,3,4,6} provides a framework for assessing the quality, reliability, and interpretability of clustering results. Cluster validation determines whether the obtained clusters depict meaningful underlying structure in the data. In some variants of Mapper, clusters are formed at fixed intervals, regardless of whether an intrinsic cluster structure exists; therefore, validation is essential for comparing clustering algorithms, selecting the number of clusters, tuning parameters, and detecting overfitting or spurious partitions.

The theoretical foundation of cluster validation is primarily based on two criteria: compactness and separation~\cite{2}. Compactness requires that data points belonging to the same cluster be close to each other according to a chosen similarity or distance measure, while separation requires that data points belonging to different clusters be well separated. A clustering solution is considered adequate when it achieves a suitable balance between these two objectives.
Internal validation indices such as the Silhouette score~\cite{1}, Dunn index~\cite{2}, and Davies-Bouldin index~\cite{3} assess this balance using only the data and the clustering result through geometric or statistical measures of cohesion and dissimilarity. External indices such as ARI~\cite{4} and NMI~\cite{43} evaluate agreement with reference labels when ground truth is available. In Mapper, these measures are computed locally according to the datasets used.

\subsection{Cluster Quality in Mapper}
\label{cluster-quality-in-mapper}
The quality of clusters produced by the Mapper algorithm depends not only on the clustering method alone but also on the other two components of the pipeline:

\begin{itemize}
    \item \textbf{Lens function:} It determines the projection of the data into a lower-dimensional parameter space on which the cover gets defined. This aligns with geometric structure, enabling meaningful separation of points within pullback sets. Therefore, a lens function affects cluster quality because an ill-suited lens may collapse distinct regions, degrading cluster quality.
    
    \item \textbf{Cover:} It identifies how the lens space is partitioned into overlapping regions, and hence it controls the resolution at which local clustering is performed. If the cover is very coarse, it can mix regions that are actually separate, and vice versa. Therefore, the cover also determines the quality of the clusters we will get.
    
    \item \textbf{Clustering algorithm:} It groups points within each pullback set. Density-based and hierarchical clustering methods make fundamentally different assumptions about how clusters should look and how they should be separated. These assumptions, in turn, interact strongly with the geometry imposed by the lens and the chosen cover.
\end{itemize}

In this direction, several works have explored fuzzy clustering based variants of Mapper, including $F$-mapper~\cite{13} and Shape Fuzzy $C$-Means (SFCM)~\cite{41}. These methods introduce soft membership assignments to improve cluster quality, but their validation criteria focus on clustering performance, with topological quality assessed by visual inspection of the resulting graphs. G-Mapper~\cite{65} adopts a statistically guided approach to cover optimization, in which cover elements are iteratively refined using normality testing and Gaussian mixture models. Although this method does not define an explicit validation index, it aligns the cover construction with the underlying data distribution. A more systematic framework is proposed in D-Mapper~\cite{35}, which bridges clustering and topology by combining the Silhouette Coefficient with a persistent homology-based measure, termed the Topological Signal Rate. This measure quantifies significant topological features using bottleneck bootstrap confidence sets. It allows Mapper outputs to be compared across different parameter settings. The results show that optimizing clustering quality in isolation is not sufficient to preserve meaningful topological structure. Overall, existing Mapper variants illustrate a clear progression from clustering focused validation toward approaches that incorporate statistical and topological considerations. The development of a cluster-quality measure focused on Mapper graphs is still lacking.

\subsection{Cluster Quality Experiments}
\label{sec:cluster_quality_experiments}

Before analyzing the results, it is important to clarify how we measure "quality" in this context. The metrics presented here are calculated on \textit{Voronoi-induced hard partitions} of the Mapper nodes, rather than on the overlapping graph itself. Since Mapper nodes share data points~\cite{10,23}, we assign each point to a single node based on the nearest centroid.

This simplification allows us to use standard validation metrics like the Silhouette Coefficient and Normalized Mutual Information (NMI)~\cite{1,4}. However, readers should be careful while interpreting the results; these metrics measure the \textit{geometric compactness} of the simplified partition, not the topological faithfulness of the Mapper graph. Consequently, a high Silhouette score might simply point out that the algorithm has fragmented a continuous shape into small, spherical clusters. This is actually a topological failure. Conversely, low scores may correctly reflect the topology of the manifold. Hence, we treat these metrics as diagnostic tools for assessing parameter sensitivity, not as objectives to be optimized.

\subsubsection*{1. Dataset Comparison}
\label{sec:dataset_comparison}

\begin{figure}[H]
  \centering
  \includegraphics[width=0.85\textwidth]{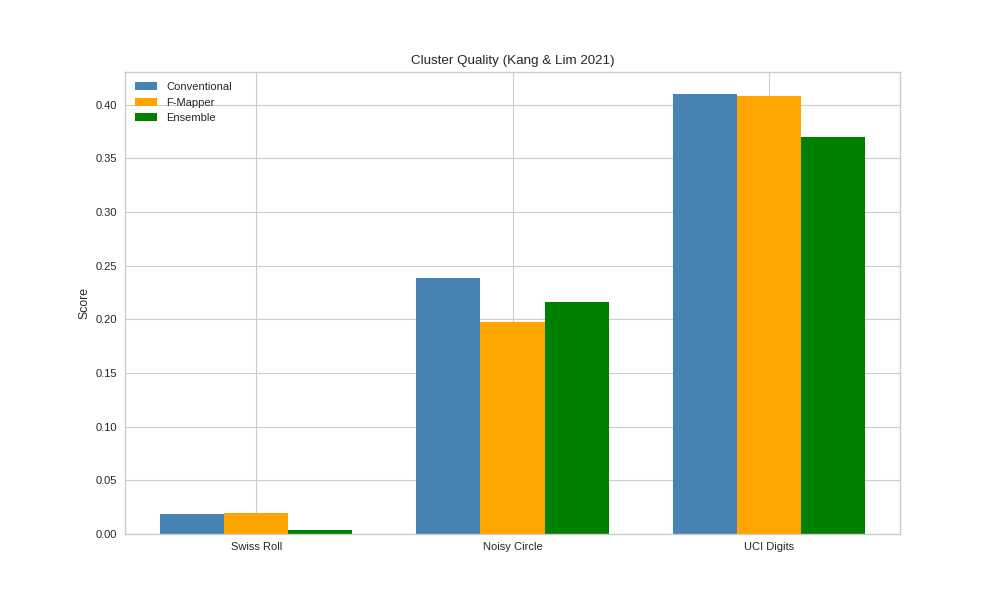}
  \caption{\textbf{Cluster quality across datasets.} Silhouette Coefficient (SC) for Swiss Roll and Noisy Circle, Normalized Mutual Information (NMI) for UCI Digits. Low SC values on the Swiss Roll reflect manifold geometry and the effect of Voronoi hardening; the Ensemble is slightly more conservative in these runs. Numerical values are reported in Table~\ref{tab:cluster_quality_results}.}
  \label{fig:dataset_comparison}
\end{figure}

We evaluated three variants of Mapper: Conventional, $F$-mapper, and Ensemble across three datasets, each chosen for its distinct topological properties (Figure~\ref{fig:dataset_comparison}).

\paragraph{Swiss Roll (Continuous Manifold).}
The continuous 2D manifold embedded in 3D yields consistently very low silhouette values for all methods tested. Measured scores are: Conventional $\approx 0.018$, $F$-mapper $\approx 0.020$, and Ensemble $\approx 0.003$. These near-zero values indicate that Voronoi hardening produces partitions that do not align with compact, spherical clusters on the manifold (as expected). The Ensemble is slightly more conservative here and produces the smallest silhouette, which we interpret as an avoidance of artificial fragmentation rather than a failure of the method.

\paragraph{Noisy Circle (Discrete Structure with Noise).}
For the Noisy Circle we observe moderate silhouette values, consistent with a geometry amenable to geometric clustering. Measured scores are: Conventional $\approx 0.238$, $F$-mapper $\approx 0.197$, and Ensemble $\approx 0.216$. These values indicate that all methods capture a meaningful circular structure to varying degrees; the Conventional Mapper in our runs shows the highest silhouette, while the Ensemble provides a compromise between robustness and geometric compactness.

\paragraph{UCI Digits (High-Dimensional Classification).}
On the UCI Digits dataset (NMI against ground-truth labels) the measured scores are: Conventional $\approx 0.410$, $F$-mapper $\approx 0.408$, and Ensemble $\approx 0.370$. Here the Ensemble is slightly more conservative and attains marginally lower NMI than single-run variants in these particular runs. Importantly, ensemble diagnostics (selected grid parameters, selected silhouette values, node counts, final $K$) are saved and available in the repository for inspection; these diagnostics explain the conservative behaviour observed in some ensemble runs.

\begin{table}[h]
\centering
\caption{Cluster Quality Comparison (Silhouette for synthetic datasets; NMI for UCI Digits). Values reflect the runs after enforcing Min--Max scaling and the ensemble fixes described in the codebase.}
\label{tab:cluster_quality_results}
\begin{tabular}{lccc}
\toprule
\textbf{Dataset} & \textbf{Conventional} & \textbf{F-mapper} & \textbf{Ensemble} \\
\midrule
Swiss Roll (Sil) & 0.018 & 0.020 & 0.003 \\
Noisy Circle (Sil) & 0.238 & 0.197 & 0.216 \\
UCI Digits (NMI) & 0.410 & 0.408 & 0.370 \\
\bottomrule
\end{tabular}
\end{table}

\subsubsection*{2. Parameter Sensitivity Analysis}
\label{sec:parameter_sensitivity}

\begin{figure}[H]
  \centering
  \includegraphics[width=1.1\textwidth]{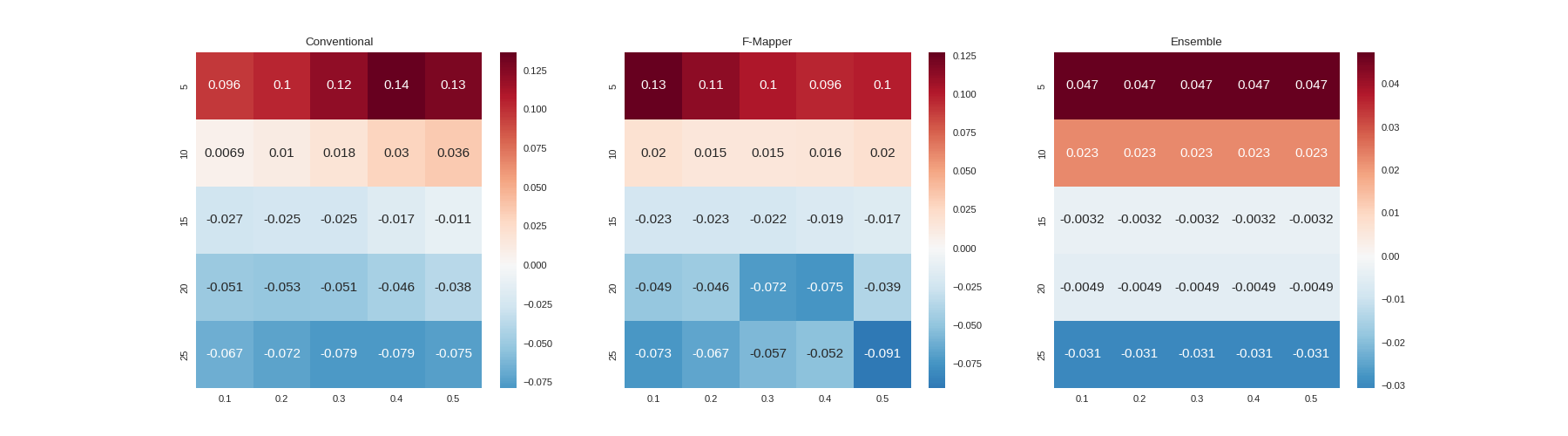}
  \caption{\textbf{Parameter sensitivity on Swiss Roll.} Heatmaps show Silhouette scores as functions of resolution $n$ and overlap $p$. Conventional and $F$-mapper methods degrade due to fragmentation at high resolution; the Ensemble reduces volatility by selecting consensus maps.}
  \label{fig:parameter_heatmaps}
\end{figure}

To understand the stability of the effect of parameter changes on Silhouette and NMI, we varied the resolution $n \in \{5, 10, \dots, 25\}$ and overlap $p \in \{0.1, \dots, 0.5\}$ on the Swiss Roll (Figure~\ref{fig:parameter_heatmaps}), keeping DBSCAN parameters and the lens function (PCA1) fixed.

\paragraph{Conventional Mapper (Figure~\ref{fig:parameter_heatmaps}).}
Performance generally degrades as resolution increases: the best scores are observed at coarser covers, and finer covers tend to reduce silhouette due to geometric fragmentation. This behaviour is consistent with the fragmentation phenomenon described in Part I.

\paragraph{$F$-mapper (Figure~\ref{fig:parameter_heatmaps}).}
$F$-mapper shows a similar degradation with increasing resolution. While fuzzy intervals can help locally, they do not prevent fragmentation at high resolutions in this manifold setting, and silhouette values remain low across much of the grid.

\paragraph{Ensemble Mapper (Figure~\ref{fig:parameter_heatmaps}).}
The Ensemble tends to produce more stable (less volatile) silhouette responses across the scanned grid. Scores remain low overall because the manifold geometry resists compact partitioning; however, the ensemble reduces dramatic parameter-driven swings by selecting consensus maps.

\subsubsection*{3. Effect of Lens Function}
\label{sec:lens_effect}

\begin{figure}
  \centering
  \includegraphics[width=0.60\textwidth]{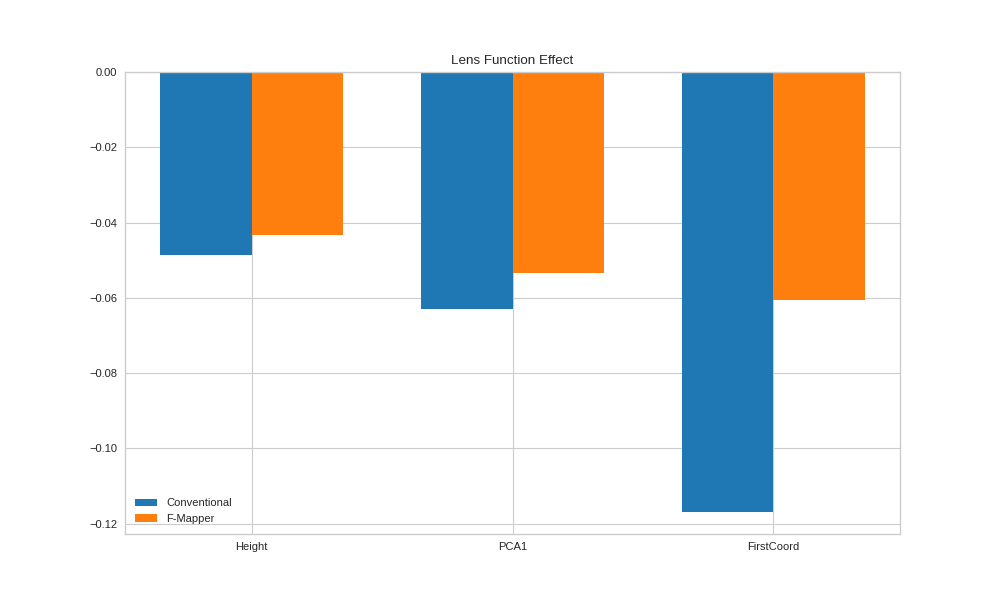}
  \caption{\textbf{Effect of lens function on Swiss Roll.} All lenses resulted in low or negative silhouette values, confirming the manifold's continuous nature and the limits of Voronoi hardening.}
  \label{fig:lens_effect}
\end{figure}

Lens selection is often one of the most important user choices~\cite{8,10}. We tested three lenses on the Swiss Roll (Figure~\ref{fig:lens_effect}), holding other parameters constant.

\paragraph{Height (Second Coordinate).}
This lens produced low or slightly negative silhouette values in our runs (Conventional: $-0.049$, $F$-mapper: $-0.043$). The height coordinate is multi-valued along the roll, which creates long, non-convex Voronoi regions and hence poor silhouette scores.

\paragraph{PCA1 (First Principal Component).}
PCA1 also yields low scores (Conventional: $-0.063$, $F$-mapper: $-0.053$), as it compresses the 3D structure into a single axis causing distant spiral segments to be projected nearby.

\paragraph{First Coordinate Projection.}
First-coordinate projection shows mixed performance (Conventional: $-0.117$, $F$-mapper: $-0.060$). Overall, discrete clustering metrics remain poor on this manifold and should be interpreted as diagnostic of geometric mismatch rather than as indicators of algorithmic failure.

\subsubsection*{4. Effect of DBSCAN Parameters}
\label{sec:dbscan_effect}

We investigated the impact of DBSCAN's radius threshold ($\varepsilon$), varying it from 0.1 to 0.7.

\begin{figure}[H]
  \centering
  \begin{subfigure}{0.48\textwidth}
    \centering
    \includegraphics[width=\linewidth]{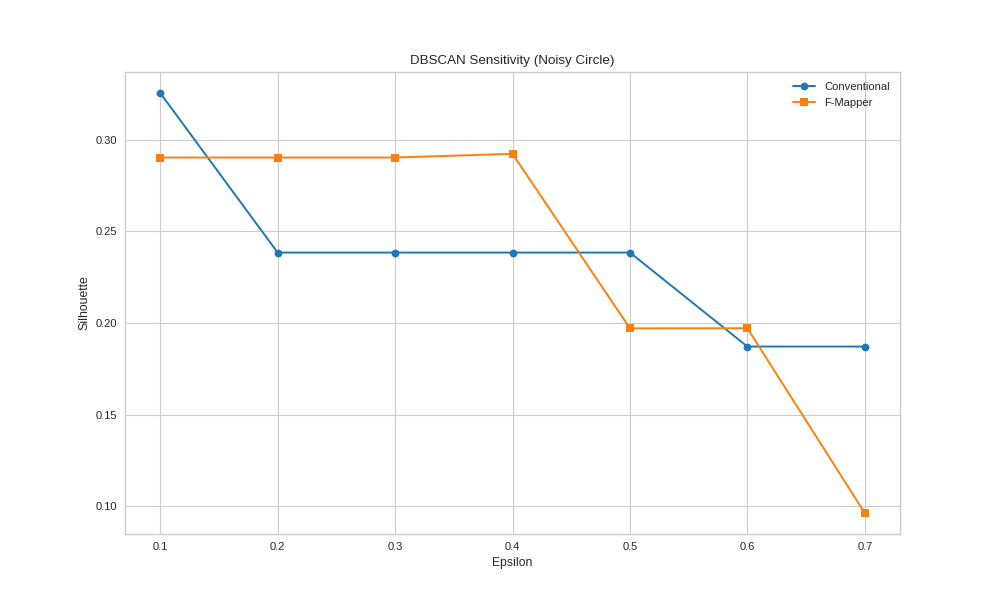}
    \caption{Noisy Circle}
    \label{fig:dbscan_circle}
  \end{subfigure}
  \hfill
  \begin{subfigure}{0.48\textwidth}
    \centering
    \includegraphics[width=\linewidth]{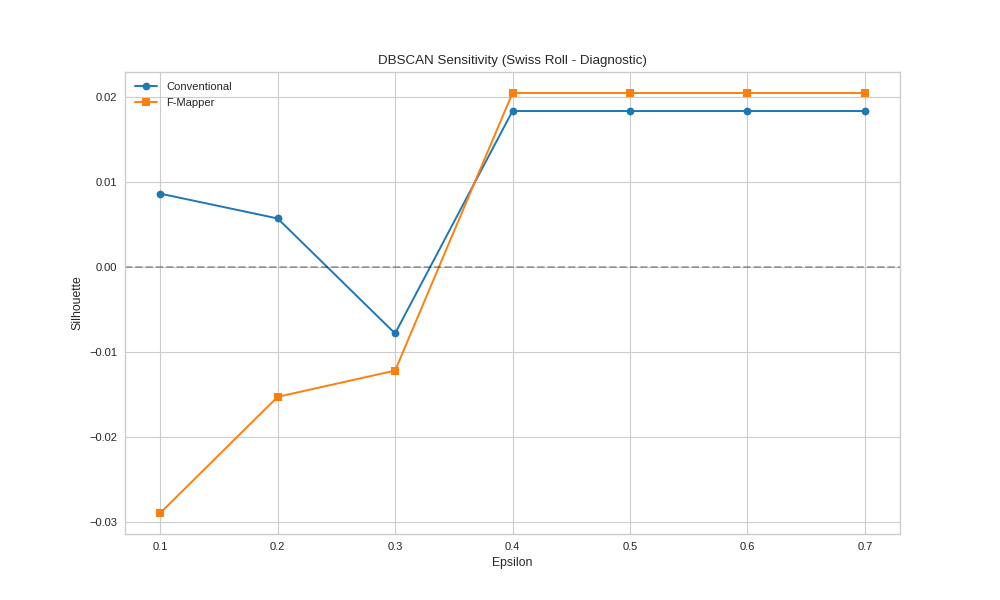}
    \caption{Swiss Roll (diagnostic sweep)}
    \label{fig:dbscan_swiss}
  \end{subfigure}
  \caption{\textbf{Effect of DBSCAN $\varepsilon$ on cluster quality.} (a) The Noisy Circle shows a stable plateau with a clear regime of good silhouette. (b) The Swiss Roll exhibits uniformly low silhouette values across parameters; small positive values at larger $\varepsilon$ do not indicate preserved topology.}
  \label{fig:dbscan_effect}
\end{figure}

\paragraph{Noisy Circle (Figure~\ref{fig:dbscan_circle}).}
For the Noisy Circle, we observe a reasonably stable regime. In our runs, the Conventional Mapper achieved a peak silhouette near $\varepsilon=0.1$ (Conventional $\approx 0.326$), then plateaued around $\approx 0.238$ for larger radii; the $F$-Mapper attains values around $0.29$ across a range of $\varepsilon$ values. These results indicate that with appropriate $\varepsilon$, geometric cluster compactness can be recovered for discrete circular structure.

\paragraph{Swiss Roll (Figure~\ref{fig:dbscan_swiss}).}
For the Swiss Roll, silhouette values remain very small across the sweep (e.g., Conventional values range from $\approx 0.009$ to $\approx 0.018$, and fuzzy/overlap-aware scores fluctuate around small negative or near-zero values). Unlike some extreme runs reported with different preprocessing, here we do not observe large spuriously-high silhouette spikes; instead the manifold remains resistant to DBSCAN-based, compact clustering except at parameter extremes. This again highlights that improved silhouette does not necessarily imply preserved topology.

\subsubsection*{Summary and Interpretation}

This analysis reflects the limits of using geometric metrics to check topological structures:

\begin{enumerate}
\item \textbf{Manifolds resist discrete clustering.} The persistently low silhouette scores on the Swiss Roll are diagnostic of continuous topology and the incompatibility with Voronoi-induced hard partitions.

\item \textbf{Context matters for parameters.} Discrete datasets (Noisy Circle) show stable parameter regimes where geometric clustering aligns well with intuition; continuous manifolds (Swiss Roll) do not. Parameter selection should therefore be informed by domain knowledge and topological goals.

\item \textbf{Ensemble conservatism.} The Ensemble method tends to produce more conservative (less extreme) geometric scores in our runs. Where the Ensemble returns slightly lower silhouette/NMI, diagnostics indicate that this behavior stems from its consensus selection and node-clustering steps; the saved JSON diagnostics can be inspected for exact selected parameters and membership matrices.

\item \textbf{Diagnostic value.} The DBSCAN experiments emphasize that apparently "better" geometric scores can correspond to over-fragmentation or loss of topological structure; these results motivate the topological shape-preservation checks presented in Part III.
\end{enumerate}

\subsection{Open Problems}
\label{open-problems}
In the above sections, we have demonstrated the concept of cluster quality. In particular, in Section~\ref{sec:cluster_quality_experiments}, we discussed how changes in different parameters affect clustering performance. Based on these observations, we highlight several open problems that remain to be addressed:

\begin{itemize}
    \item Developing a framework that jointly optimizes all relevant parameters, rather than tuning them independently.
    \item Identifying appropriate validation metrics that provide a more comprehensive and reliable characterization of the underlying datasets.
    \item Most cluster quality metrics used in Mapper are borrowed from classical clustering and are primarily geometric in nature; therefore, there is a need for cluster quality measures that assess topological structure rather than purely geometric properties.
\end{itemize}

\section*{
\begin{center}
PART III 
\end{center}
}

\section{Topological and Shape Preservation}
\label{sec:topological_shape_preservation}

The term \emph{preservation} in Mapper-based methods is used in multiple senses, often conflated. In this part, we explicitly distinguish between two fundamentally different notions:
\begin{itemize}
    \item \textbf{Topological Preservation}, understood relative to Reeb graphs or Reeb spaces induced by a function.
    \item \textbf{Shape Preservation}, understood as stability of metric and connectivity structure independent of any function. It is the stability of the $\varepsilon$-net nerve, not the recovery of homotopy type.
\end{itemize}

This distinction is essential to avoid conceptual ambiguity. Conventional Mapper and its theoretical extensions aim to approximate Reeb-theoretic objects associated with a chosen lens, whereas Ball Mapper is a purely metric construction designed to preserve geometric shape and connectivity rather than function-induced topology. This section is further divided into several subsections.
Section~\ref{sec:reeb_perspective} discusses topology preservation from a Reeb graph perspective;
Section~\ref{sec:ball_mapper_shape} examines shape preservation;
Section~\ref{sec:topology_experiments} presents the related experiments;
and Section~\ref{sec:preservation-open-problems} discusses the corresponding open problems.

\subsection{Reeb Graphs Perspective for Ground Truth}
\label{sec:reeb_perspective}

Given a topological space $X$ and a continuous function $f: X \to \mathbb{R}$, the Reeb graph encodes the evolution of connected components of level sets of $f$. Reeb graphs arise naturally in Morse theory~\cite{66} and provide a canonical low-dimensional descriptor of function-induced topology. For multivariate functions $f: X \to \mathbb{R}^d$, the Reeb graph generalizes to the \emph{Reeb space}, which captures the connectivity of joint level sets. Although Reeb spaces are difficult to compute directly, they constitute the natural continuous objects that Mapper constructions aim to approximate~\cite{9}. Accordingly, \emph{topological preservation for Mapper is mainly assessed relative to the Reeb graph or Reeb space induced by the same lens function}.

\paragraph{Topology Preservation via Reeb Graph}

The structural connection between Mapper and Reeb graphs has been rigorously established for one-dimensional functions. In particular, the MultiNerve Mapper can be interpreted as the Reeb graph of a perturbed function~\cite{26}.
A key consequence of this interpretation is that Mapper does not introduce spurious topological features at the continuous level. Instead, it selectively retains or discards features of the Reeb graph depending on the interaction between:
\begin{itemize}
    \item the placement and resolution of the cover in the codomain, and
    \item the location of critical values of the function.
\end{itemize}
Loss or distortion of topological features is therefore predictable and governed by discretization parameters rather than arbitrary algorithmic instability. Consequently, Carrière and Oudot~\cite{26} formally demonstrate that Mapper behaves as a pixelized or discretized version of the Reeb space.

\paragraph{Structural Relationship Between Mapper and Reeb Graphs}

Beyond structural correspondence, Mapper has been shown to be a statistically consistent estimator of the Reeb graph under suitable assumptions~\cite{8}. As the sample size increases and the cover resolution is refined in a controlled manner, Mapper converges to the Reeb graph in bottleneck distance via extended persistence representations~\cite{8}.
Probabilistic convergence results further establish that Mapper approximates the Reeb graph with high probability when data are sampled from noisy distributions. In these settings, approximation error is explicitly bounded by the cover resolution, linking discretization choices to topological accuracy~\cite{29}.
These results clarify that deviations from topological preservation observed in finite samples reflect resolution limits rather than violations of theoretical guarantees.

\paragraph{Categorical Convergence to Reeb Spaces}
Further work is provided by categorical approaches, where both Mapper and Reeb spaces are analyzed as functors encoding connected components over open sets~\cite{9}.
In this framework, Mapper converges to the Reeb space in interleaving distance, with error bounded by the resolution of the cover. This applies to multivariate functions and unifies previous structural and statistical analyses.
In the scalar case, categorical convergence implies geometric convergence, reinforcing the interpretation of Mapper as a principled discretization of Reeb-theoretic objects.

\paragraph{Limitations of Reeb-Theoretic Preservation}

Despite strong theoretical guarantees, Reeb-based preservation is inherently conditional on the choice of function. Topological features not aligned with the selected lens may be suppressed, fragmented, or absent in the resulting Mapper graph.
Moreover, at high resolutions, discretization may induce combinatorial artifacts such as cycle proliferation. These artifacts do not correspond to genuine Reeb-space features and motivate consideration of alternative notions of preservation that do not depend on a filter function~\cite{23,26}.

\subsection{Shape Preservation via Ball Mapper}
\label{sec:ball_mapper_shape}

Ball Mapper~\cite{12} provides a complementary notion of preservation based purely on the metric structure of the data and not the topology. It constructs a graph as the nerve of an $\varepsilon$-net cover of the point cloud, where vertices correspond to balls and edges indicate shared coverage.
Unlike Mapper, Ball Mapper does \emph{not} approximate Reeb graphs or Reeb spaces. It therefore does not aim to recover intrinsic Betti numbers of the underlying manifold. Instead, it preserves:
\begin{itemize}
    \item local neighborhood geometry,
    \item global connectivity,
    \item and stability of the covering structure across scales.
\end{itemize}
Critically, Ball Mapper's cycle counts primarily reflect covering density and nerve connectivity rather than intrinsic topological features. For this reason, Ball Mapper is most appropriately interpreted as a \emph{geometric stability baseline} rather than a topological ground truth.

\paragraph{Relationship Between Mapper and Ball Mapper}

When the lens function used in Mapper is uniformly continuous, partial correspondences between Mapper and Ball Mapper graphs can be established. Points that are close in Ball Mapper are expected to map to nearby regions in Mapper~\cite{12}. However, the converse does not generally hold. Ball Mapper may preserve metric connectivity that Mapper intentionally separates due to function-induced structure. Consequently, Ball Mapper may obscure functional topology while providing a more faithful representation of geometric shape. Ball Mapper should therefore be interpreted not as a replacement for Mapper, but as a complementary geometric baseline emphasizing metric stability over Reeb theoretic fidelity. Finally, Dłotko~\cite{12} establishes a relation between Ball Mapper and the Vietoris-Rips complex.

\subsection{Topology-Shape Preservation Experiment}
\label{sec:topology_experiments}

In this section, we concentrate our experimental analysis on \textbf{Conventional Mapper}, \textbf{$F$-mapper}, and \textbf{Ball Mapper}. These three mappers correspond to the different paradigms of standard Reeb approximation, fuzzy clustering, and geometric covering, respectively.
We eliminate \textbf{Multiscale Mapper} and \textbf{Ensemble Mapper} from this comparison because of methodological problems:
\begin{enumerate}
    \item \textbf{Multiscale Mapper} outputs a higher-order simplicial complex rather than a simple 1-skeleton graph. Applying the graph-theoretic cycle rank formula ($\beta_1 = E - V + C$) used in this section would yield mathematically invalid comparisons, as it neglects the cycle-collapsing effect of higher-dimensional simplices (see Table~\ref{tab:variant_evaluation_axes}).
    \item \textbf{Ensemble Mapper} is a meta-algorithmic architecture that aims to counteract instability. Its behavior is better understood in the context of \textit{consensus stability}, which is the specific concern of Part IV (Section~\ref{sec:integrative}).
\end{enumerate}
  
We have analyzed experiments on two bases:
\begin{itemize}
    \item \textbf{Topology Preservation:} The ability to recover the ground-truth Betti numbers of the Reeb graph (skeleton).
    \item \textbf{Shape Preservation:} The ability to capture the geometric surface of the underlying manifold~\cite{12,28}.
\end{itemize}

\subsubsection*{1. Mapper vs Reeb Graph (Ground Truth Topology)}

We evaluated whether Mapper-based constructions approximate the Reeb graph induced by a chosen lens function (1D height projection). \textbf{For this specific lens}, the Reeb graph has a simple 1-dimensional skeleton with $\beta_1 = 1$ (the hole; the extrinsic spiral collapses under height projection)~\cite{26}. This provides a \textit{lens-specific ground truth} distinct from the Euclidean persistent homology ($\beta_1 \in \{1,2\}$) computed in Part~I.

\textbf{Critical Note on Betti Number Calculation:}
In Part~III, the value of $\beta_1$ is calculated in the \textbf{1-skeleton} of the Mapper graph with \textit{graph-theoretic cycle rank} (Equation: $b_1 = E-V + C$). Note that in contrast to the simplicial homology computation, in which filled-in triangles are treated as null homologous, every possible cycle in the wireframe is counted. Within the framework of ``Topological Explosion''~\cite{23}, every value obtained serves a critical purpose.

Table~\ref{tab:topology_comparison} is a comparison table for Conventional Mapper and $F$-mapper. Conventional Mapper has a large degree of instability, with the value for $\beta_1$ averaging out to $16.7$ ($1567\%$ deviation, ground truth $= 1$). $F$-mapper reduces this value to $4.5$ ($350\%$ deviation), yet both appear unstable.
\begin{table}[!ht]
\centering
\begin{tabular}{cccc}
\toprule
\textbf{Resolution ($N$)} & \textbf{Conventional $\beta_1$} & \textbf{$F$-mapper $\beta_1$} & \textbf{Status} \\
\midrule
6  & 1  & 1 & Stable \\
12 & 1  & 3 & Stable \\
24 & 26 & 11 & Diverging \\
28 & 28 & 6 & Exploded (Conv) \\
32 & 26 & 4 & Exploded (Conv) \\
40 & 18 & 2 & Exploded (Conv) \\
\bottomrule
\end{tabular}
\caption{\textbf{Swiss Roll dataset: Comparison of $\beta_1$ values across methods.}
With the increase in resolution, the $\beta_1$ parameter of the Conventional Mapper starts from $1$ and suddenly increases to $28$, whereas the $F$-mapper shows more robustness in this increase. Analyzing the figures more closely, it can be realized that the volatilities of the Conventional Mapper are higher compared to the $F$-mapper because the standard deviation of $\beta_1$ for the Conventional Mapper is $\sigma_{\beta_1} = 11.05$, whereas that of the $F$-mapper is $\sigma_{\beta_1} = 2.93$. (Data Source: \texttt{part3\_quantitative\_metrics.csv}).}
\label{tab:topology_comparison}

\end{table}

\subsubsection*{2. When High Resolution Triggers Topological Explosion}

Among the most problematic phenomena we've observed in the Conventional Mapper, certainly one of the most interesting is the phenomenon we've termed ``Topological Explosion.'' Here is the way in which it occurs: We start increasing the resolution of the filter, and as a consequence, the graph, rather than being a topological representation, becomes something much more like a geometric mesh~\cite{23}.

\begin{figure}[!ht]
    \centering
    \includegraphics[width=\linewidth]{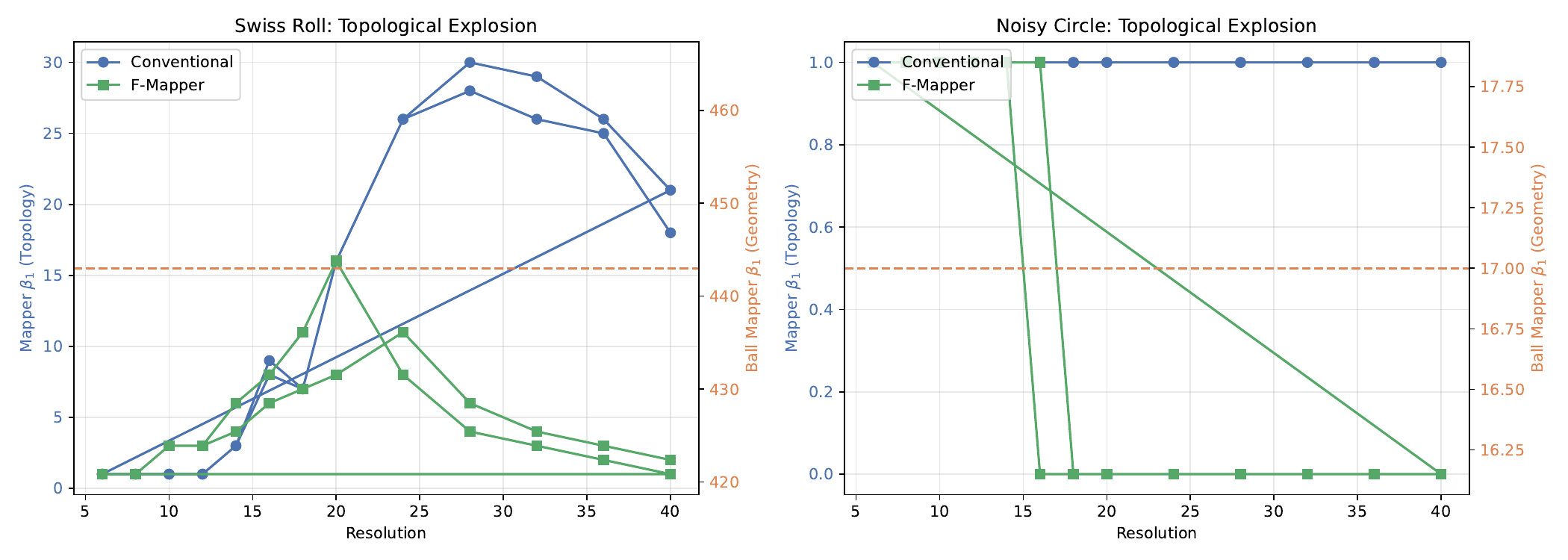}
    \caption{\textbf{How resolution drives topological explosion.} 
    The blue curve (Conventional Mapper) shows spurious loops multiplying as resolution increases. $F$-mapper (green) reduces this problem substantially. Ball Mapper (orange dashed line) maintains a constant but very high Betti number ($\beta_1 \approx 455$), because it is capturing the geometric lattice covering the surface rather than adapting to resolution changes in the skeletal structure. (Data Source: \texttt{part3\_topological\_preservation.csv}).}
    \label{fig:topo_explosion}
\end{figure}

In Figure~\ref{fig:topo_explosion}, one can observe this phenomenon. The blue plot shows the number of $\beta_1$ in the Conventional Mapper steadily increasing with increasing resolution. Rather than decoding the topological features, this algorithm appears to compute the fuzzy tessellation of the geometry. Then, the green plot indicates $F$-mapper, which uses fuzzy clustering~\cite{13,41} to mitigate this explosion problem. On the other hand, Ball Mapper, represented by the orange dashed line, is opposite, as its $\beta_1$ is constant independently of resolution ($\sigma_{\beta_1}=0.0$). However, this constancy is obtained at the cost of meaningfulness, as the number that is kept constant ($\beta_1 \approx 455$) represents the lattice structure on the surface, not the topology.

\paragraph{How Different Methods Handle Noisy Data.} 
When we examined how these algorithms respond to noise, we found an interesting contrast. Conventional Mapper turns out to be quite resilient to small amounts of noise; the change in $\beta_1$ is only about $0.77$ when we increase noise from $\sigma=0.02$ to $0.05$. This robustness comes from how the method works: its grid-based structure effectively rounds nearby points into the same bins. Ball Mapper tells a different story entirely. it is much more sensitive to noise (with $\beta_1$ changes around $12.0$ for the same noise increase), which makes sense given that the method is designed to ensure complete coverage of all data points, including outliers~\cite{12}.

\subsubsection*{3. Testing How Filter Choice Affects Results (Topology vs Shape)}

The differences we see between these methods really come down to a fundamental tension: should we prioritize recovering Reeb graph structure, or should we focus on capturing geometric coverage?

\begin{itemize}
    \item \textbf{Conventional Mapper} relies heavily on the chosen filter function~\cite{26}. When we increase the resolution, the way intervals overlap rigidly can fragment the data, turning the output into what looks like a ``mesh'' rather than a clean topological summary.
    
    \item \textbf{$F$-mapper} uses fuzzy smoothing to become less sensitive to filter choices~\cite{13}, but this creates its own problem: on datasets with simple topology (like the Noisy Circle), it can smooth things out too much. We've seen cases where $\beta_1$ drops to zero at high resolution values of $N$.
\end{itemize}

\subsubsection*{4. Comparing Geometric Coverage: Ball Mapper vs Filter-Based}

\begin{figure}[!ht]
    \centering
    \includegraphics[width=\linewidth]{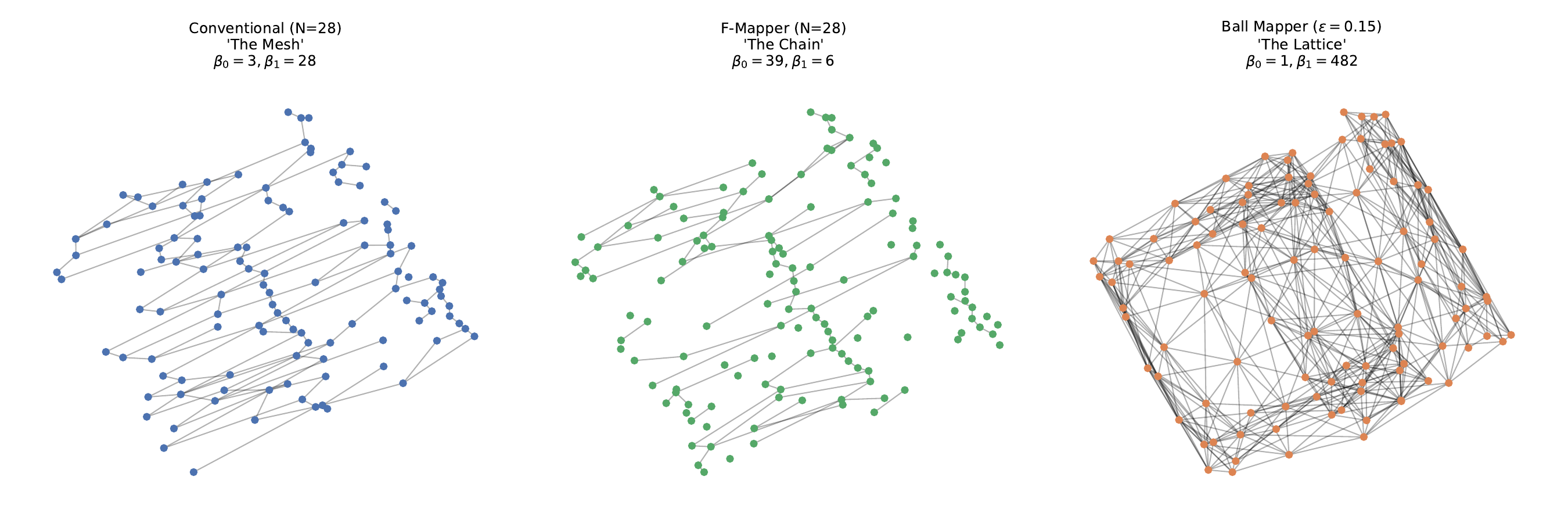}
    \caption{\textbf{Comparing graph structures at $N=28$.} 
    On the left, Conventional Mapper creates a mesh-like structure ($E/V \approx 1.2$). In the middle, $F$-Mapper distills things down to a skeletal chain. On the right, Ball Mapper wraps the manifold surface in a dense lattice ($E/V \approx 8.1$, with $\beta_1 \approx 482$).}
    \label{fig:structure_comparison}
\end{figure}

This is an important point to note that Conventional Mapper and $F$-mapper are obviously trying to reconstruct the topology~\cite{26}, Ball Mapper has a different objective altogether, which is to try and place a geometric skin around the data~\cite{12}. We can see the difference in Figure~\ref{fig:structure_comparison}, for example, in the output produced for $N = 28$. To make the comparison more rigorous, we computed the following measure for the graph density, which we define as the ratio of edges $E$ to vertices $V$.

\begin{itemize}
    \item \textbf{Conventional ($E/V \approx 1.2$) and $F$-mapper ($E/V \approx 0.9$):} These produce sparse, tree-like graphs; they're building skeletal structures.
    
    \item \textbf{Ball Mapper ($E/V \approx 8.1$):} This creates a much denser, lattice-style graph, wrapping a skin around the data.
\end{itemize}

The fact that the density in Ball Mapper is about 8 times higher confirms what we suspected, it is essentially a geometric covering algorithm. Therefore, when we notice it is producing large Betti numbers, it is obviously not a bug; it is the algorithm doing what it is supposed to.

\paragraph{The Connectivity-Cycle Tradeoff.} 
There is another interesting pattern here: these methods behave in opposite ways when it comes to $\beta_0$ (the number of connected components). Conventional Mapper only fragments the data slightly ($\beta_0$ averages around $2.3$), but $F$-Mapper fragments it aggressively as a strategy for maintaining topological stability.
\begin{itemize}
    \item \textbf{$F$-mapper's typical behavior:} Averages about $\beta_0 \approx 20.6$ connected components across different resolutions.
    
    \item \textbf{$F$-mapper at $N=28$:} As we can see in Figure~\ref{fig:structure_comparison}, the fragmentation becomes extreme, reaching $\beta_0 = 39$ separate components.
\end{itemize}

Ball Mapper, by contrast, keeps everything connected in a single component ($\beta_0 = 1.0$). This tells us something important about how $F$-mapper achieves its stability: it suppresses spurious loops ($\beta_1$) by chopping the data into many small, isolated clusters. This also explains something we noticed back in Part~II: $F$-mapper tends to get high Silhouette scores because it creates tight, well-separated clusters. The downside? It sacrifices global connectivity to achieve this~\cite{4,13}.

\subsubsection*{5. Validation of Theorem 3.4 (Interleaving)}

Our work aimed at applicable validation of Theorem~3.4 proposed by D\l{}otko et al.~\cite{12}. This statement has a theoretical meaning, as it claims that the graph $G(\varepsilon)$ obtained via the Ball Mapper should be woven among the persistence modules of the Vietoris-Rips complex $VR(r)$.

\begin{figure}[!ht]
    \centering
    \begin{subfigure}{\linewidth}
        \centering
        \includegraphics[width=1.0\linewidth]{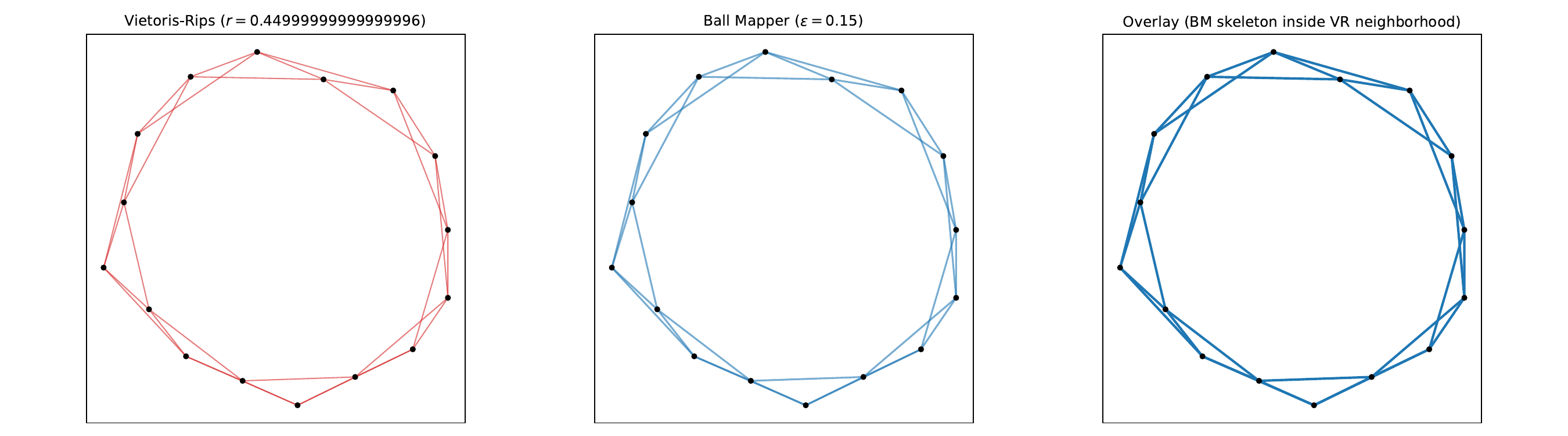}
        \caption{Noisy Circle (2D): Ball Mapper (Blue) forms a skeleton inside the Vietoris-Rips complex (Red).}
    \end{subfigure}
    \par\bigskip 
    \begin{subfigure}{\linewidth}
        \centering
        \includegraphics[width=1.0\linewidth]{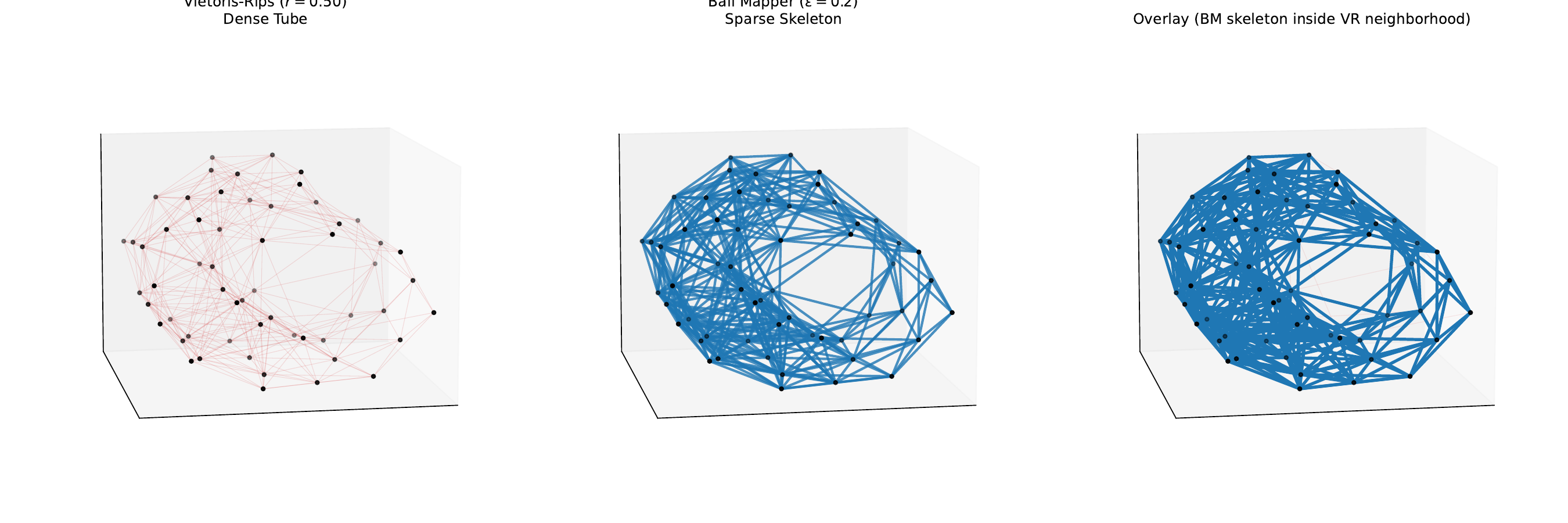}
        \caption{Swiss Roll (3D): Ball Mapper lattice (Blue) tracks the Rips manifold (Red).}
    \end{subfigure}
    \caption{\textbf{Visual Intuition for Theorem~3.4 (Interleaving).} The blue Ball Mapper graph is geometrically contained within the red Vietoris-Rips complex, consistent with theoretical interleaving bounds~\cite{12}.}
    \label{fig:theorem_3_4}
\end{figure}

This structural relationship is illustrated in Figure~\ref{fig:theorem_3_4}. This figure is indicative of the existence of a very important distinction that can be proved via experiments: even though the high values of the Betti numbers in the case of Ball Mapper are caused by the surface lattice effect, there is no arbitrariness in the cycles detected in the manifold, and the technique is very efficient when the scale $r$ is larger than $2\varepsilon$~\cite{10,12}.

\subsection{Open Problems in Topological and Shape Preservation}
\label{sec:preservation-open-problems}

Although substantial progress has been made in understanding the preservation properties of Mapper and its variants, several important open problems remain.
\begin{itemize}
    \item Existing theoretical results establish asymptotic convergence of Mapper to Reeb graphs or Reeb spaces under idealized conditions. However, explicit and computable bounds on discretization error for finite samples, especially under noisy and irregular sampling, remain largely unexplored.
\item Mapper and Ball Mapper represent fundamentally different preservation objectives, i.e., Reeb-theoretic topology versus metric shape. Developing hybrid or multi-objective frameworks that balance topological correctness with geometric coverage remains an open research direction.
\item Betti numbers alone are often insufficient to characterize meaningful preservation, particularly in the presence of dense combinatorial artifacts. More expressive invariants and summaries that distinguish genuine topological structure from discretization-induced noise are needed.
\item 
    The robustness of topological and shape preservation under varying noise levels differs significantly across Mapper variants. A systematic theory explaining how noise interacts with cover resolution, clustering, and nerve construction is still missing.
\item While Reeb-theoretic interpretations are well developed for scalar lenses, preservation guarantees for multivariate lenses and Reeb spaces are far less mature, especially in practical, high-dimensional settings.
\item  There is currently no unified metric that jointly evaluates topology preservation relative to Reeb graphs and shape preservation in the metric sense. Designing such metrics would enable principled comparison across Mapper variants and parameter regimes.
\end{itemize}
These open problems emphasize that preservation in Mapper-based methods is inherently multi-faceted, and that reliable interpretation requires joint consideration of topology, geometry, and discretization effects.

\section{Integrative Comparison and Applications}
\label{sec:integrative}
In this section we analyze the relation between these parts; Section \ref{sec:stability_cluster_quality} depicts the relation between \ref{Part I - Stability Perspective}
and section \ref{Part II - Cluster Quality Perspective}, section \ref{Relation of Cluster Quality with Topology Preservation} depicts relation between section \ref{Part II - Cluster Quality Perspective} and section \ref{sec:topological_shape_preservation}, section \ref{Relation of Preservation with Stability} analyzes the relation between section \ref{sec:topological_shape_preservation} and section \ref{Part I - Stability Perspective} and section \ref{sec:stability_ensemble} gives a brief analysis of stability of Ensemble Mapper as promised above. Section \ref{Decision Framework for Method Selection} gives an idea for selection of Mapper variants and section \ref{sec:integrative-open-problems}  summarizes the open problems identified throughout the paper.

\subsection{Relation of Stability with Cluster Quality}
\label{sec:stability_cluster_quality}

This subsection establishes the connection between the stability properties of Mapper constructions \ref{Part I - Stability Perspective} and the cluster quality measures analyzed in Part~II.

\begin{figure}
    \centering
    \includegraphics[width=0.4\linewidth]{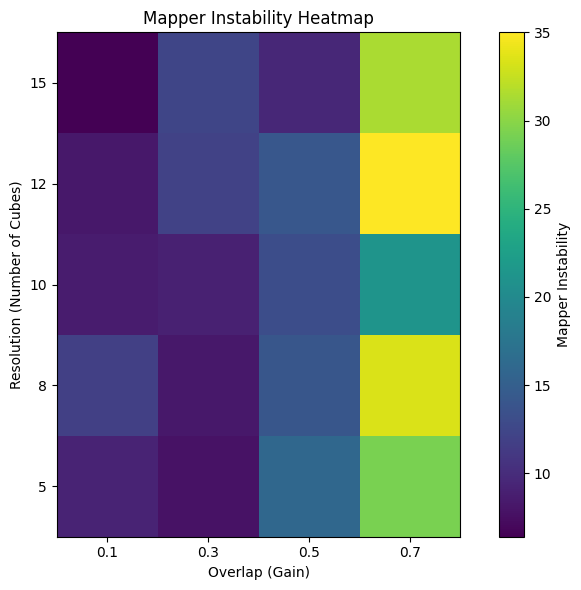}
   \caption{\textbf{Conventional Mapper Instability Heatmap.}
This heatmap is obtained for the Swiss Roll dataset using the same parameter
settings as in Fig.~\ref{fig:parameter_heatmaps}.}

    \label{Fig 12}
\end{figure}

Intuitively, one might expect a stable Mapper construction to produce clusters whose quality does not degrade significantly under small perturbations. However, our experimental results contradict this intuition and indicate that stability and cluster quality are not fully dependent on each other. From Fig.~\ref{Fig 12}, we observe that the parameter setting $n=15, p=0.1$ is relatively stable compared to other configurations, yet Fig.~\ref{fig:parameter_heatmaps} shows that the corresponding silhouette score is not the highest. Coversely, high silhouette score does not imply low instability. This demonstrates that stability is neither a necessary nor a sufficient condition for improved clustering performance. Consequently, stability and cluster quality should be treated as separate axes when analyzing datasets, rather than as complementary or interchangeable.

\subsection{Relation of Cluster Quality with Topology Preservation}
\label{Relation of Cluster Quality with Topology Preservation}

\begin{figure}[h]
    \centering
    \includegraphics[width=0.4\linewidth]{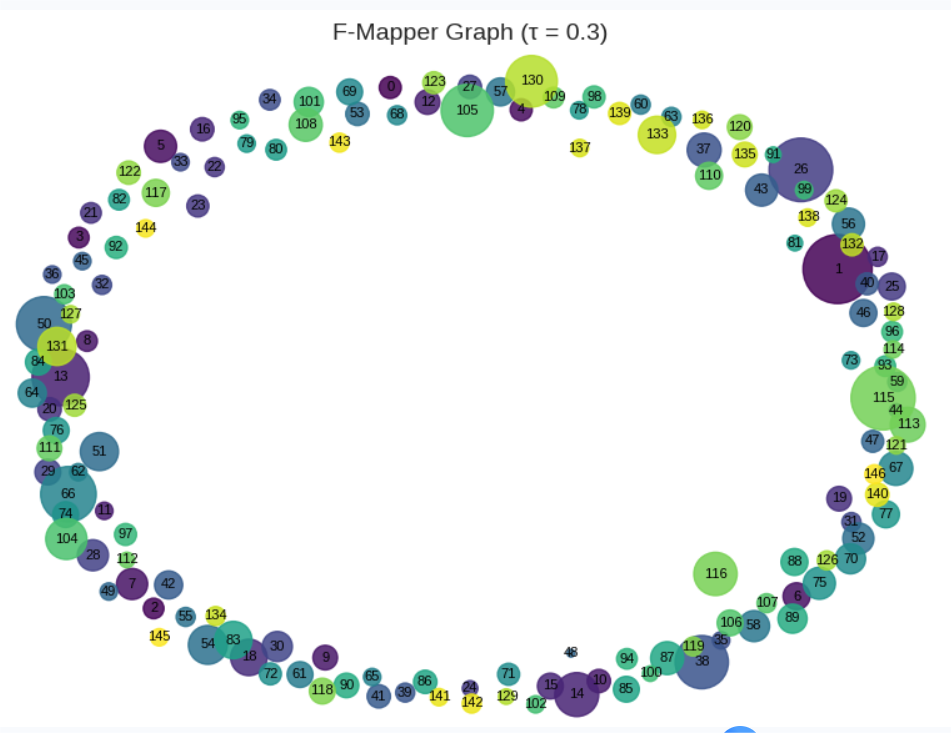}
   \caption{\textbf{$F$-mapper on the Swiss Roll Dataset.}
This result is obtained using the parameter settings as in with $n = 10$ and $\tau = 0.3$.}

    \label{$F$-mapper on Swiss Roll Dataset}
\end{figure}

Figure~\ref{$F$-mapper on Swiss Roll Dataset} is obtained by varying the number of intervals 
$n \in \{5, 10, 15, 20, 25\}$ and the overlap percentage 
$\tau \in \{0.05, 0.10, \ldots, 0.30\}$, and then selecting the parameter $\tau$ based on the silhouette score. 
The figure clearly demonstrates poor topological reconstruction due to disconnected components. 
This observation further supports the claim that a high silhouette score does not necessarily imply 
preservation of the underlying topology or geometric shape of the data.

\subsection{Relation of Preservation with Stability}

\label{Relation of Preservation with Stability}
\begin{itemize}
    \item 
If we look at Ball Mapper's Bootstrap Stability in Fig \ref{fig:bootstrap_stability}, we find that Ball Mapper's coefficient of variation is 8.8 \%, which is comparatively lower than Conventional and $F$-mapper.  But in Fig \ref{fig:geometric_heatmap}, it is showing that $\beta_1$ varies from 339-482. On the Swiss Roll, Ball Mapper consistently produces hundreds of 1-cycles. This stable artifact reflects the packing density of the balls on the manifold surface. While geometrically consistent (the lattice persists across noise levels), it is
topologically distinct from the underlying manifold.
\item In figure \ref{Fig 12}, the parameter setting $n$ = 12, $p$ = 0.7 yields the best stability score, but in Fig \ref{17}, we can easily visualize that neither the topology nor the shape is preserved. This is due to the effects of the parameters we use throughout the process.

\begin{figure}
    \centering
    \includegraphics[width=0.5\linewidth]{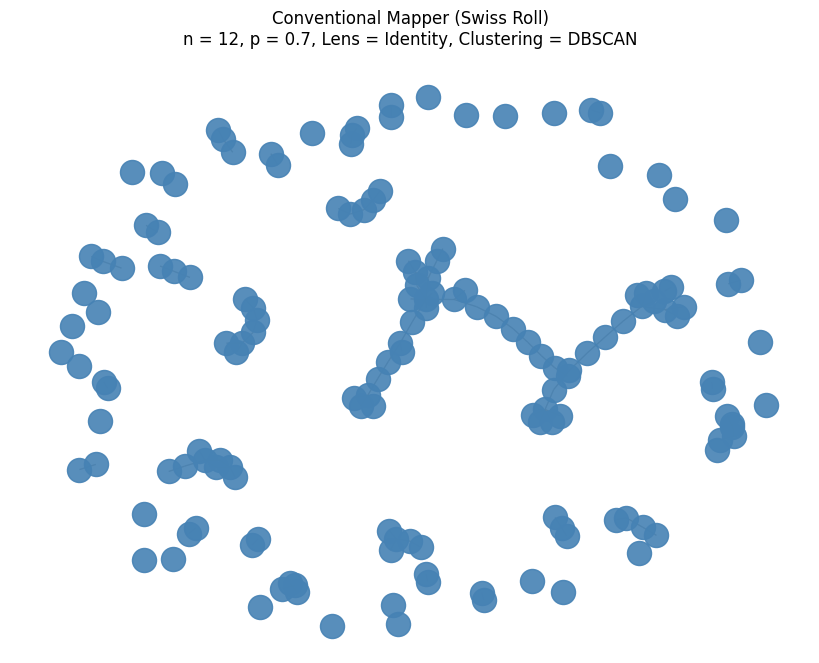}
    \caption{\textbf{Conventional Mapper Visualization of the Swiss Roll Dataset.}
The graph is obtained using the First coordinate projection lens function $f(x,y)=x$,
with $(n,p)=(12,0.7)$ and DBSCAN as the pullback clustering algorithm.}
\label{17}
\end{figure}
\end{itemize}

\subsection{Stability Analysis of Ensemble Mapper}
\label{sec:stability_ensemble}

As outlined in the preliminaries, this section presents the stability analysis of the Ensemble Mapper. To evaluate the stabilizing capability of ensemble approaches, we implemented the \textbf{Ensemble Mapper (Algorithm 2)} proposed by Kang and Lim~\cite{11}. Unlike simple averaging strategies, this method uses a \textbf{Selection-and-Consensus} approach:
\begin{enumerate}
    \item \textbf{Grid Search:}  For a target resolution $N$, we generate a pool of base partitions by varying the interval count $L \in [N-2, N+2]$ and overlap $p \in [20\%, 40\%]$.
    \item \textbf{Selection:} We select the $k=10$ highest-quality partitions using the Silhouette Coefficient~\cite{1}. This step filters out poorly structured (e.g., shattered) partitions from the ensemble.
    \item \textbf{Meta-Clustering:} A similarity matrix is generated between each Mapper node from the selected partitions. These are then subjected to hierarchical clustering to identify stable "meta nodes," ultimately leading to the creation of the Ensemble Graph.
\end{enumerate}

We compared this approach against a standard single-run Conventional Mapper across our three datasets.

\begin{figure}[!ht]
    \centering
    \begin{subfigure}{0.48\textwidth}
        \includegraphics[width=\linewidth]{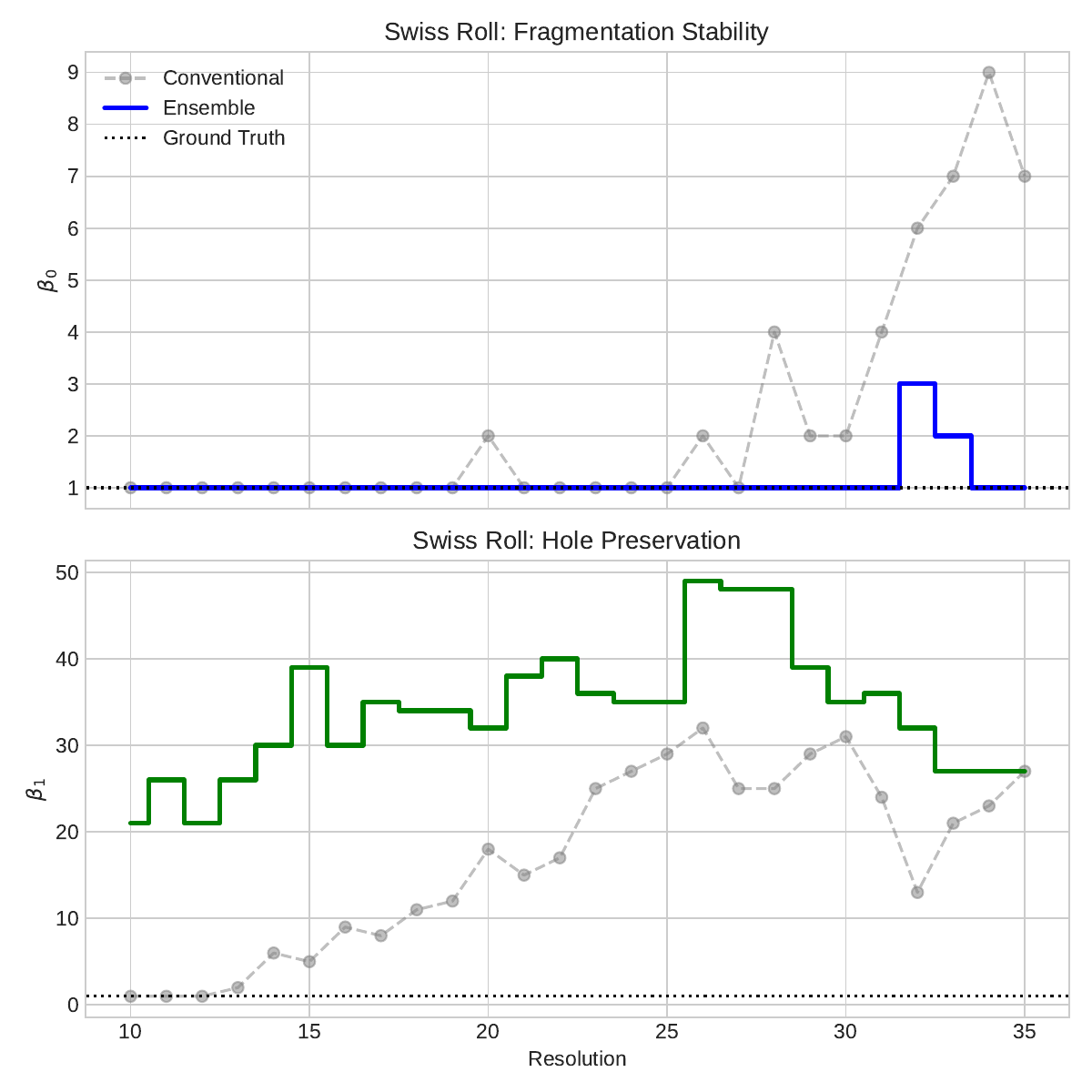}
        \caption{Swiss Roll (Continuous)}
        \label{fig:part4_swiss}
    \end{subfigure}
    \hfill
    \begin{subfigure}{0.48\textwidth}
        \includegraphics[width=\linewidth]{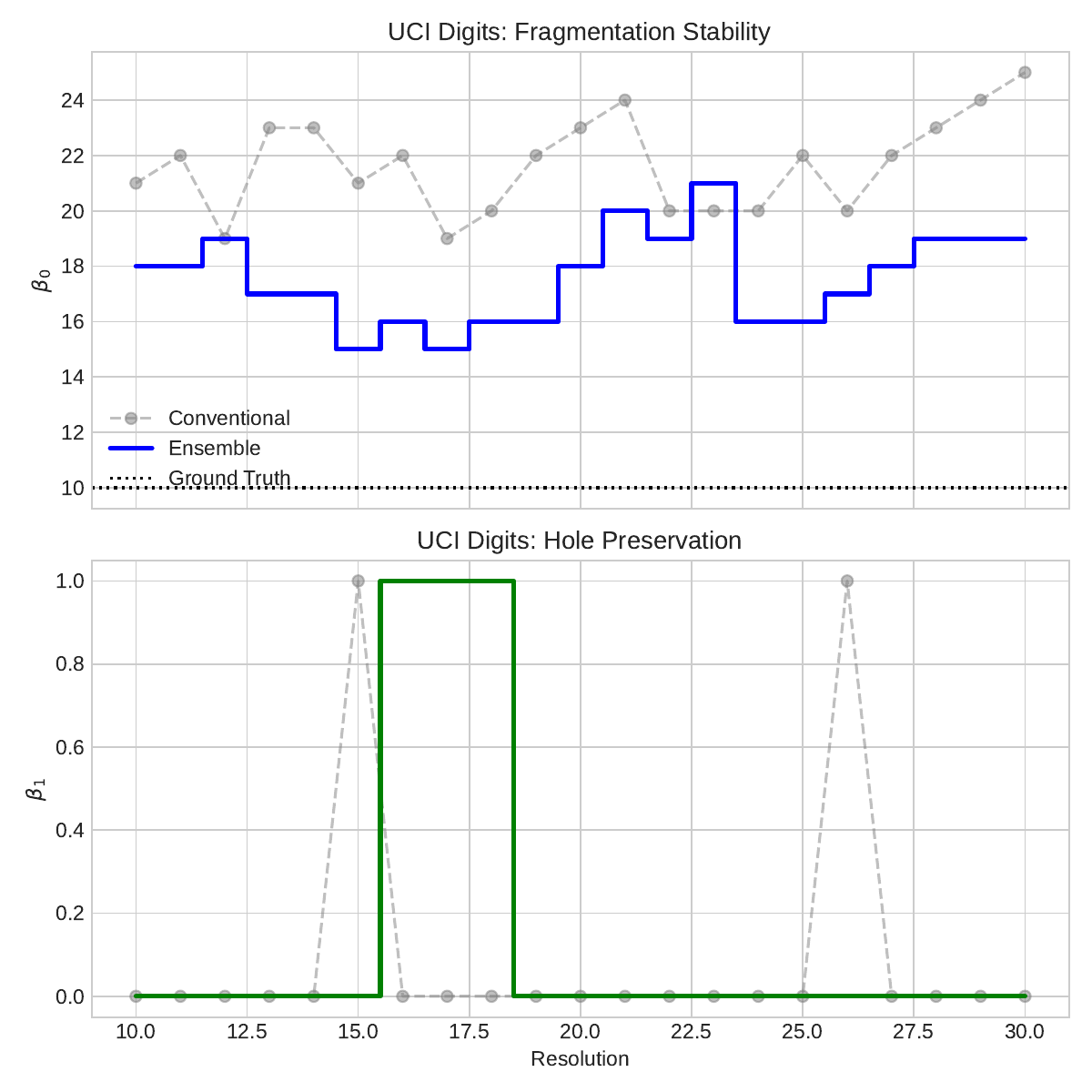}
        \caption{UCI Digits (High-D)}
        \label{fig:part4_digits}
    \end{subfigure}
    
    \caption{\textbf{Integrative Stability Analysis.} 
  (a) The single-run Mapper (grey) on the Swiss Roll exhibits a topological explosion ($\beta_0$ rises to 9) at high resolution. The Ensemble (blue) prevents this explosion, maintaining a constant $\beta_0 = 1$, although the $\beta_1$ count is higher due to the "thick" meta-graph structure. 
(b) On UCI Digits, the Ensemble handles fragmentation caused by sparsity better than the single-run Mapper, maintaining lower $\beta_0$ and greater stability.}
    \label{fig:ensemble_stability}
\end{figure}

The results (Figure~\ref{fig:ensemble_stability} and Table~\ref{tab:stability_results}) reveal a distinct trade-off in the ensemble approach:

\begin{enumerate}
    \item \textbf{Mitigation of Fragmentation ($\beta_0$):}
    A significant advantage of the Ensemble method is its ability to mitigate topological explosion~\cite{23,26}. On the Swiss Roll, the Conventional Mapper creates numerous disconnected components (nine) at the highest resolution, while the Ensemble method maintains a single component across all dimensions (evidenced by $\beta_0=1$; the standard deviation was reduced from 2.33 to 0.43).
Similarly, on the UCI Digits dataset, where increased sparsity causes disconnection in the single-run Mapper, the Ensemble maintains stable connectivity and creates compact structures resembling the original classes.
    
    \item \textbf{The ``Mesh Effect'' on Homology ($\beta_1$):}
    While the Ensemble preserves the \textit{existence} of 1-dimensional holes (never dropping to 0), it tends to overestimate their number. On the Swiss Roll, $\beta_1$ stabilized at high values ($\approx 27$) rather than the ground truth of 1. On the simple Noisy Circle, the Ensemble unnecessarily complicated the graph ($\beta_1 > 40$). This artifact arises because the meta-clustering process creates a ``thick,'' triangulated graph (a mesh) rather than a clean skeleton.
\end{enumerate}

\begin{table}[h]
\centering
\caption{Stability Comparison (Standard Deviation of Betti Numbers across Resolution Sweep). Lower indicates greater stability.}
\label{tab:stability_results}
\begin{tabular}{lcccc}
\toprule
& \multicolumn{2}{c}{\textbf{Fragmentation ($\sigma_{\beta_0}$)}} & \multicolumn{2}{c}{\textbf{Hole Stability ($\sigma_{\beta_1}$)}} \\
\textbf{Dataset} & \textbf{Conventional} & \textbf{Ensemble} & \textbf{Conventional} & \textbf{Ensemble} \\
\midrule
Swiss Roll & 2.33 & \textbf{0.43} & 10.37 & 7.44 \\
Noisy Circle & \textbf{0.00} & 0.00 & \textbf{0.00} & 15.74 \\
UCI Digits & 1.71 & \textbf{1.66} & 0.36 & 0.36 \\
\bottomrule
\end{tabular}
\end{table}

These findings support a nuanced conclusion: \textit{The Kang \& Lim Ensemble acts as a powerful topological regularizer for connectivity ($\beta_0$), effectively curing topological explosion in complex manifolds. However, for 1-dimensional features ($\beta_1$), it trades precision for persistence, often producing robust but topologically complex meta-graphs.}

\subsection{Decision Framework for Method Selection}

\label{Decision Framework for Method Selection}
This subsection provides a brief guide on when to use each Mapper variant, and Table~\ref{tab:mapper_comparison} provides a comprehensive comparison of Mapper variants across the three evaluation axes examined in this review.
\begin{enumerate}
    
\item \textbf{When to use Conventional Mapper:}
\begin{itemize}
    \item Exploratory analysis with interpretable lens functions.
    \item Fast prototyping and visualization.
    \item When complemented with stability and cluster quality validation.
\end{itemize}

\item \textbf{When to use MultiNerve Mapper:}
\begin{itemize}
    \item Scientific discovery requiring theoretical guarantees.
    \item When correct Betti numbers are essential.
    \item Publications requiring rigorous stability claims.
\end{itemize}

\item \textbf{When to use Fuzzy Mapper:}
\begin{itemize}
    \item Noisy data with uncertain boundaries.
    \item Anisotropic or nonlinear lens structures (use GK or kernel variants).
    \item When smooth transitions are preferred over hard partitioning.
\end{itemize}

\item \textbf{When to use Ball Mapper:}
\begin{itemize}
    \item Lens-free analysis of intrinsic geometry.
    \item When no natural filter function exists.
    \item Topology recovery on manifolds.
\end{itemize}

\end{enumerate}

\begin{table}[H]
\centering
\caption{Comprehensive comparison of Mapper variants across evaluation axes}
\label{tab:mapper_comparison}
\small
\begin{tabular}{@{}llll@{}}
\toprule
\textbf{Variant} & \textbf{Strengths} & \textbf{Weaknesses} & \textbf{Best Use} \\
\midrule
Conventional & Simple; fast; & Parameter sensitive; & Exploratory \\
Mapper & lens-interpretable & spurious cycles & with validation \\
\midrule
MultiNerve & Theoretical stability; & Computational & Scientific \\
Mapper & correct Betti numbers & overhead & discovery \\
\midrule
Fuzzy & Reduces boundary & Higher quality $\neq$ & Noisy \\
Mapper & artifacts; smooth & better topology & boundaries \\
\midrule
Ball & Preserves geometry; & Radius selection & Lens-free \\
Mapper & lens-independent & critical & analysis \\
\midrule
Ensemble & Averages instability; & Multiple runs; & Exploratory \\
Mapper & identifies robust features & consensus may mask & work \\
\bottomrule
\end{tabular}
\end{table}

\subsection{Open Problems}
\label{sec:integrative-open-problems}
Despite recent advances in the theoretical and empirical analysis of Mapper and its variants,
several fundamental open problems remain unresolved when the three evaluation axes-stability,
cluster quality, and topological or shape preservation are considered jointly.
\begin{itemize}
\item Mapper parameters are tuned
independently with respect to stability, clustering quality, or topology preservation.
Developing a principled framework that jointly optimizes all three axes or characterizes
their Pareto trade-offs remains an open challenge.

Most existing analyses evaluate Mapper outputs at the graph level, ignoring variability within
individual cover elements. Our experiments indicate that instability, poor clustering, and
topological explosion often originate locally at the interval level. Formal diagnostics that quantify interval-wise failure modes are still lacking.

\item The choice of lens function strongly influences stability, cluster compactness, and topological fidelity, yet no theory explains how to select lenses that balance these competing objectives. Characterizing lens functions that are jointly stable and topology-preserving is an open
problem.

\item  Topological explosion is often treated as an unavoidable consequence of increasing cover resolution. Designing Mapper variants or regularization strategies that suppress the proliferation of combinatorial cycles without sacrificing resolution or interpretability remains unresolved.

\item 
While $F$-mapper and Ensemble Mapper empirically mitigate instability and topological explosion, their theoretical behavior across the three axes is poorly understood. Establishing convergence, stability bounds, and failure modes for these methods is an open area of research.

\item Ball Mapper exhibits strong geometric stability but does not aim to recover Reeb-theoretic topology, while Conventional Mapper targets Reeb graphs but is unstable at high resolutions. Resolving this gap by designing methods that preserve both metric shape and functional topology remains open.

\item Our results show that no single Mapper variant performs optimally across all datasets and evaluation axes. Developing data-driven selection criteria for Mapper variants is an important open direction.

\item Despite heavy use, Mapper evaluation still relies on visual inspection.
Designing quantitative, interpretable metrics that simultaneously reflect stability, cluster quality, and topological faithfulness remains a major open challenge.

\end{itemize}

\section{Conclusion}
\label{sec:conclusion}

Mapper has proved to be one of the most popular tools for analyzing topological data, thanks to its graph-based approach for reducing data dimensionality. Although it is widely applied, our discussion of the Mapper process across various design parameters remains incomplete and fragmented in the literature. The review presents an attempt to organize these concepts into three aspects: stability, cluster quality, and topology and shape preservation.

In Part I, we focused on stability from both theoretical and experimental perspectives. We demonstrated how stability cannot be assessed by looking at how well the answer remains comparable for varying parameters. Indeed, a Mapper construction may be relatively stable yet fail to reveal interesting topological information. In Part~II, we investigated cluster quality in Mapper and demonstrated how sensitive it is to parameters such as the cover parameter settings, lens functions, and the clustering algorithm. We also demonstrated how cluster validity indices are not sound indicators of how well the topology is being reconstructed. In Part~III, we investigated topological and shape reconstruction and emphasized the distinction between Reeb topological reconstruction and shape reconstruction.

In Part IV, we have combined these viewpoints and examined stability, cluster quality, and preservation simultaneously. Our results show that they frequently conflict with each other. It might happen that optimizing with regard to one criterion leads to deteriorated performance with respect to another. For instance, having a good silhouette value does not necessarily imply good topological preservation; it might happen that a stable Mapper output has too much or misleading topological complexity, or that a geometric shape-preserving function might fail to detect topological structure.

In conclusion, this review makes it clear that Mapper should not be viewed as a rigid algorithm in itself, but rather as a flexible system that, depending on a number of interrelated design decisions, can produce different outcomes. In order to properly interpret the data yielded by a Mapper analysis, it is essential to be familiar with interdependencies among the techniques. We believe that with the open issues raised throughout this work, research techniques will be developed in order to create more meaningful applications of Mapper, which will finally be applicable in real-world scenarios.

\section*{Disclosure statement}\label{sec:disclosure-statement}
No potential conflict of interest was reported by the author(s).

\section*{Author contributions}\label{sec:author-contributions}
Annesha Sen and Shivam Singh contributed equally to the conceptualization of the problem and the execution of the computations. S.P. Tiwari supervised the work and helped in the preparation of manuscript.

\section*{Statements of ethical approval}
For this research article, the authors did not undertake work that involved human participants or animals.

\end{document}